\documentclass[conference]{IEEEtran}
\IEEEoverridecommandlockouts
\usepackage{cite}
\usepackage{amsmath,amssymb,amsfonts}
\usepackage{algorithmic}
\usepackage{graphicx}
\usepackage{textcomp}
\usepackage{xcolor}
\usepackage{bm}
\usepackage{epstopdf}
\usepackage{lineno,hyperref,float,amsfonts,amsmath,epsfig,graphics,booktabs,multirow,cite,stfloats}
\usepackage{comment}
\def\BibTeX{{\rm B\kern-.05em{\sc i\kern-.025em b}\kern-.08em
    T\kern-.1667em\lower.7ex\hbox{E}\kern-.125emX}}
\begin{document}

\title{Active learning for anti-disturbance dual control of unknown nonlinear systems
	
	\thanks{This work was supported by the National Natural Science Foundation of China (Grant no. 62073259).}
}
\author{Xuehui Ma, Shiliang Zhang, Fucai Qian, Jinbao Wang, Yushuai Li
	
	\thanks{Xuehui Ma is with the School of Automation and Information Engineering, Xi'an University of Technology, Xi'an, China (e-mail: maxuehuiphd@stu.xaut.edu.cn). }
	\thanks{Shiliang Zhang is with the Department of Informatics, University of Oslo, Oslo, Norway (e-mail: shilianz@ifi.uio.no).}
	\thanks{Fucai Qian is with the School of Automation and Information Engineering, Xi'an University of Technology, Xi'an, China (e-mail: qianfc@xaut.edu.cn).}
	\thanks{Jinbao Wang is with the School of Automation and Information Engineering, Xi'an University of Technology, Xi'an, China (e-mail: jbwang@stu.xaut.edu.cn). }
	\thanks{Yushuai Li is with the Department of Informatics, University of Oslo, Oslo, Norway (e-mail: yushuaili@ieee.org).}
}

\newtheorem{assumption}{Assumption}[section]
\newtheorem{lemma}{Lemma}[section]
\newtheorem{remark}{Remark}[section]
\newtheorem{thm}{\bf Theorem}[section]

\maketitle

\begin{abstract}
This work concerns the control of unknown nonlinear systems corrupted by disturbances. For such systems, we propose an anti-disturbance dual control approach with active learning of the disturbances. Our approach holds the dual property of handling the two tasks simultaneously and iteratively: (i) learn the disturbances affecting the system and (ii) drive the system output towards a reference trajectory. Particularly, we model nonlinear system dynamics using a specialized neural network (SNN). This SNN formulates the disturbances via the designed additive and multiplicative disturbance components. We consider both additive and multiplicative disturbances for precise description and recognition of disturbance profile. We achieve the disturbance recognition in the SNN via the design of a Bayesian-based active learning approach, which allows the disturbance learning to be decoupled from the control law derivation. Such a decoupling contributes to the control robustness in the existence of varying and abrupt disturbances. We derive the dual control law based on the active learning of the SNN, and validate our approach via one-time and Monte Carlo simulations. The results demonstrate a fast disturbance recognition by our method in real-time and the robustness of control of unknown systems with abrupt disturbances. We evaluate our approach on the speed control of high-speed train, and the results manifest efficient control of the train speed with disturbance resilience, without prior knowledge about the train dynamics and the disturbances imposed on the train. We openly released the code for this work for reproduction purpose\footnote{\url{https://github.com/MXH0113/active-learning-control}}.
\end{abstract}


\section{Introduction} \label{section1}

Disturbances exist in most practical control processes in various formats, ranging from interrupting factors from the external environment to uncertainties in the inner system \cite{chen2015disturbance,8666991,8796434}. E.g., in the speed control for high-speed trains, there can be disturbances from the external environment (wind, rain, snow, etc.) and inner components (the wear, aging, and malfunctions of mechanical and electronic components of the train) \cite{5871322,6099627,8941083}. Disturbances adversely affect the system control process, which can lead to the risk of productivity reduction and economic loss. As such, suppressing the negative influence of disturbances and improving the system control robustness become crucial objectives in the control law design, especially for practical systems~\cite{DBLP:journals/tac/MaZLQSH24,DBLP:conf/eucc/MaCZLQS24,yu2024adaptive,MA2022110}. It is well known that feed-forward control strategies could eliminate the influence of measurable disturbances \cite{guzman2021tuning}. However, quite often, disturbances in practice are randomly occurred and prohibitive to be measured directly. Such a challenge motivates the development of various robust and adaptive control methods to attenuate disturbance's impact on the control~\cite{guo2014anti,li2014disturbance}.

Amongst potential solutions to the disturbance challenge, robust control strategies derive conservative control law based on the worst-case control performance in the existence of disturbances \cite{abu2006nonlinear}. For instance, Zhang et al. developed an $H_2 / H_{\infty} $ model predictive control for a partially unknown Markov jump system with bounded disturbance\cite{ZHANG20183423}. Hooshmandi et al. proposed a robust $H_{\infty}$ control method based on the polynomial linear parameter varying model for nonlinear sample-data systems \cite{LPV2018}. Reinforcement learning is applied to the $H_{\infty}$ control of unknown systems with disturbances \cite{kiumarsi2017h}. Liu et al. designed an adaptive neural network prescribed performance bounded $H_{\infty}$ controller for stochastic nonlinear systems\cite{8793219}. The robust control strategies assume a bounded scale for the disturbance and calculate a maximum cost index for the control performance under the bounded disturbance. By minimizing the maximum cost index, one obtains the robust control law. However, the obtained control law can be over-conservative due to worst-case-oriented computation \cite{qian2010adaptive,huang2016robust}.

Unlike robust control, adaptive control strategies try to learn the system dynamics and disturbance and use the learned information to adjust the control law, leading to a relatively aggressive control process for unknown systems~\cite{gao2015neural,jana2017composite,li2002variance}. Nevertheless, the learning and controlling (aka exploration and caution) in adaptive control are conflicting~\cite{li2008optimal,1137553}. I.e., the learning side needs large excitations to the system to acquire information about how the system behaves and reacts, yet large excitations degrade the control performance. While in the other way around, low excitations could stabilize the control but restrain the opportunity to explore the system for improved control performance further. First introduced by Feldaum, dual control appears as an adaptive control strategy focusing on the mentioned conflict between learning and controlling~\cite{feldbaum1960dual}. Dual control derives control laws taking into account disturbances such as noise from the surrounding environment and model parameter fluctuations, and tries to balance learning and controlling unknown systems in an iterative way~\cite{mesbah2018stochastic,filatov2004adaptive,9189668}. In recent years, dual control theory has been further developed and successfully applied to various fields. E.g., Tutsoy et al. developed a model-free dual control scheme based on Q-learning~\cite{10055961} to balance the learning and controlling of unknown systems~\cite{9514408}. Li et al. applied dual control to a concurrent learning framework for autonomous source search in an unknown environment~\cite{li2021concurrent,chen2021dual}. Klenske et al. employed dual control theory to solve the learn-control trade-off in Bayesian reinforcement learning~\cite{klenske2016dual}; Ma et al. extended the dual control strategy to control unknown systems corrupted by severe noise and outliers~\cite{ma2022adaptive,DBLP:journals/cssc/ZhangCYZH19,DBLP:journals/tnn/ZhangCYZH18}. Liu et al. devised a fault-tolerant dual control method for stochastic systems with loss-of-control effectiveness~\cite{liu2022dual,liu2021reliable}. Dual control was also adopted to develop model predictive control (MPC) approaches~\cite{https://doi.org/10.1155/2017/9402684} for systems with unknown parameters or structures~\cite{heirung2017dual,arcari2020dual,mesbah2018stochastic}. While dual control shows promising features of handling disturbances in controlling unknown systems, however, its learning and controlling procedures are coupled~\cite{milito1982innovations}. I.e., the control law at an instant $k$ is derived based on the system dynamics and disturbances learned at the previous instant $k-1$ and vice versa. An optimal solution for the control in this case generally depends on dynamic programming, which is highly computation and memory consuming and prohibitive to implement in practical situations~\cite{fabri1998dual}. As a result, the issue of coupled learning and controlling limits dual control's application where the disturbance takes effect over time and can vary frequently, sometimes even abruptly.

To narrow this gap, constructing a structure or procedure for dual control is crucial in which the entailed disturbance learning and output controlling are decoupled. As far as we know, such a structure or procedure is not available yet. Some preliminary works have explored different ways to formulate disturbances in unknown systems, which brings deep insights into the disturbance side and how the disturbance can be possibly handled for better control. Among such efforts, most of them formulate the additive disturbance in unknown systems~\cite{8666991,8961962,sun2018robust}. The additive disturbances refer to those entering the model of a system only additively. Several studies also formulate the multiplicative disturbance to fit practical scenarios, e.g., petroleum seismic exploration~\cite{liu2015optimal}, underwater communication systems~\cite{dai2019adaptive}, radar target tracking~\cite{frost1982model}, and cyber-physical systems~\cite{song2020distributed}. In recent years, several works have already leveraged the concept of multiplicative disturbance for improved control robustness. E.g., Zhao et al. developed $H_{\infty}$ control for uncertain Markov jump systems with multiplicative noises using the linear matrix inequation technique~\cite{zhao2017robust}. Mazouchi et al. proposed a convex optimization approach using the system-level synthesis for linear quadratic regulating (LQR) with multiplicative disturbances~\cite{mazouchi2022data}. Wang et al. built a backstepping-based control for strict-feedback nonlinear systems corrupted by multiplicative noises~\cite{wang2020novel}.

This paper integrates both additive and multiplicative disturbances in the unknown system control process. We specialize a neural network to model unknown system dynamics, where the disturbances are formulated by the additive and multiplicative concepts. We embed the specialized neural network (SNN) into the dual control structure to learn and control unknown systems. This enables a derived control law with decoupled learning and controlling. Particularly, we design a strategy that assigns a collection of candidate values for the disturbance parameters in the neural network. Our strategy learns in real-time the Bayesian posterior probability for those candidates independently and in parallel. While each disturbance candidate entails a unique control law, we derive the final control law as the summation of all the entailed control laws weighted by their posterior probability. With this effort, we anticipate a control approach for unknown systems that is robust against disturbances, which take effect, vary over time, and even vary abruptly. We summarize our contribution as follows:

(i) We construct a specialized neural network (SNN) for unknown system control. We establish the SNN as an affine nonlinear neural network, and the SNN integrates both the additive and multiplicative disturbances to formulate the disturbances imposed on the system. The constructed SNN aims to precisely reflect the effect of the disturbances on the system.

(ii) We design real-time active learning of the disturbance parameters in the SNN. For this learning, we assign a finite set of candidate values to approximate the bounded disturbances, thus turning the disturbance learning into the task of recognizing probability for the possible candidates. Then, we calculate in real-time the Bayesian posterior probability for each candidate iteratively, and the disturbance parameters are recognized as the candidates' posterior probabilities converge.

(iii) Based on the designed SNN, we derive a dual control approach that decouples the learning and controlling for unknown systems. Particularly, our approach separates the disturbance learning from the system output tracking control that, the disturbance learning works in parallel with the controlling, which is a one-step-ahead optimization of trajectory tracking error and system output estimation error.

(iv) We extensively validate our approach in one-time and Monte Carlo simulations. The simulations test our approach with unknown systems corrupted by disturbances varying over time, and we impose abrupt disturbances on the systems to check the robustness of our approach. For more convincing results, we evaluate our method in the speed control of high-speed train. In those simulations and evaluations, we inspect both resilience of the control process against disturbances and how efficient the learning of disturbance parameters is.

This paper is organized as follows. Section~\ref{section2} provides the problem statement for integrating additive and multiplicative disturbances in unknown nonlinear system control. Section~\ref{section3} elaborates our anti-disturbance control strategy, including the specialized neural network (SNN), the active learning for the SNN, and the dual control law derivation for unknown systems based on the SNN. Section~\ref{section4} demonstrates simulations to validate the designed approach, evaluates our approach in the speed control of high-speed train, and compares our approach with the ideal benchmark optimal control and a latest model-free adaptive control method. We conclude our work in Section~\ref{section5}.

\section{Problem statement}\label{section2}
This section generalizes a model for unknown nonlinear systems when considering both additive and multiplicative disturbances. We explain our hypothesis for this model to be used in control procedures, and analyze the challenges of unknown system control using such a model.  

This work generalizes the objective to be controlled as a discrete-time affine nonlinear system with additive and multiplicative disturbances, as shown in the model below:
\begin{equation} \label{system}
	\begin{aligned}
		y(k+1) = \alpha(k)f[x(k)]+\beta(k)g[x(k)]u(k)+\gamma(k)+e(k),\\
	\end{aligned}
\end{equation}
where $x(k) = [y(k), y(k-1),..., y(k-n+1), u(k-1), u(k-2), ..., u(k-m)]^T$ is the system state vector, $y(k)$ is the system output, $u(k) \in R^m$ is the control signal, $f[x(k)]$ and $g[x(k)]$ are the nonlinear functions of the system state vector, $\alpha(k)$ and $\beta(k)$ are multiplicative disturbances, $\gamma(k)$ is the additive disturbance, $e(k)$ is the white Gaussian noise. We assume the following conditions for this model to represent the nonlinear system dynamics and the varying disturbances:

\begin{assumption} \label{assumption1}
	The multiplicative disturbances $\alpha(k)$, $\beta(k)$, and additive disturbance $\gamma(k)$ are unknown and time-varying in the bounded intervals $[\alpha_{l}, \alpha_{u}]$, $[\beta_{l}, \beta_{u}]$, and $[\gamma_{l}, \gamma_{u}]$, respectively. $\alpha_{l}$ and $\alpha_{u}$ are the lower and upper bound of the disturbance $\alpha(k)$, respectively. $\beta_{l}$ and $\beta_{u}$ are the lower and upper bound of $\beta(k)$, and $\gamma_{l}$ and $\gamma_{u}$ are the lower and upper bound of $\gamma(k)$, respectively.
\end{assumption}

\begin{assumption}\label{assumption2}
	The dimensionality parameters $n$ and $m$ in system state vector $x(k)$ are known.
\end{assumption}

\begin{assumption}\label{assumption3}
	The system random noise $e(k)$ follows Gaussian distribution with zero mean and variance $\sigma^2$, where the variance $\sigma^2$ is known.
\end{assumption}

\begin{assumption}\label{assumption4}
	The system is a minimum phase system, and the nonlinear function $g[ \cdot ]$ is bounded away from zero. The mathematical models of nonlinear functions $f[\cdot]$ and $g[\cdot]$ are unknown.
\end{assumption}

Based on the model and the hypothesis, we can conduct the control for unknown nonlinear systems corrupted by disturbances. For a control task of system output tracking, we formulate the dual control performance index as:
\begin{equation}\label{eq_performaceindex}
	\begin{aligned}
		J(k+1) = E\left\lbrace \sum_{k=1}^{N}[y(k+1)-y_r(k+1)]^2 |  \mathfrak{I}^k   \right\rbrace ,
	\end{aligned}
\end{equation}
where $y_r$ is the system output reference trajectory, $\mathfrak{I}^k$ is the state information at iteration $k$ and is defined as:
\begin{equation}
	\begin{aligned}
		\mathfrak{I}^k = \left\lbrace u(0), u(1), ..., u(k-1), y(1), y(2), ..., y(k) \right\rbrace .
	\end{aligned}
\end{equation}

The objective of this paper is to find a feedback control policy
\begin{equation}
	\begin{aligned}
		u(k) = \mu_k(\mathfrak{I}^k)
	\end{aligned}
\end{equation}
that minimizes the performance index $J(k+1)$ in equation (\ref{eq_performaceindex}), where the policy  $\mu(\cdot)$ is a nonlinear function to be derived to satisfy the performance index minimization.

The derivation of $\mu(\cdot)$ in real-time is challenging due to the existence of unknown disturbances denoted by $\alpha(k)$ and $\beta(k)$ and the unknown nonlinearity of $f[\cdot]$ and $g[\cdot]$. There are solutions to learn those unknown properties from the system measurements over the control process, so that an explicit control law can be derived, e.g., following the neural networks based adaptive control strategy~\cite{aastrom2013adaptive}. However, in that way, the disturbance learning and the output control are coupled. I.e., the disturbances learning at the $k$-th iteration is influenced by the control law derived at the $(k-1)$-th iteration, and the control law derivation at the $k$-th iteration depends on the learned disturbances at the $(k-1)$-th iteration. Such coupling can degrade the control of unknown systems affected by disturbances varying (even abruptly) over time. 

This paper aims to address this coupling designing a dual control strategy that breaks the coupling loop and trades off the disturbance learning and output tracking for unknown nonlinear systems represented by the model in equation~(\ref{system}). We detail our solution in the following section.

\section{The active learning for anti-disturbance dual control}\label{section3}

This section elaborates our solution of an anti-disturbance dual control for unknown nonlinear systems. We specialize a neural network to represent an unknown nonlinear system corrupted by both multiplicative and additive disturbances. We use the specialized neural network (SNN) as a basis to recognize system dynamics and disturbances that contribute to the control law derivation for unknown systems, as detailed in Section~\ref{section3.1}. Section~\ref{section3.2} illustrates how to approximate the bounded disturbances in the SNN using a finite set of candidate values, which lays the foundation for our active learning of the disturbances for unknown system control. Section~\ref{section3.3} elaborates our solution that decouples the disturbance learning and output control for unknown nonlinear systems. We demonstrate in our solution the disturbance recognition in real-time via the weight summation of the disturbance candidates and the derivation of the dual control law based on the learned SNN parameters.

\subsection{The specialized neural network for unknown nonlinear systems with disturbances}\label{section3.1}

This section specializes a neural network that complies with the model in equation~(\ref{system}). The specialized neural network (SNN) enables the numerical representation of unknown systems corrupted by disturbances, laying the foundation for system dynamics and disturbance recognition and control law derivation.

The SNN uses two Gaussian radial functions $\hat{f}$ and $\hat{g}$ to approximate the unknown nonlinear functions $f[x(k)]$ and $g[x(k)]$, respectively, where $x(k)$ is the state vector as $x(k) = [y(k), y(k-1),..., y(k-n+1), u(k-1), u(k-2), ..., u(k-m)]^T$. The two Gaussian radial functions $\hat{f}$ and $\hat{g}$ are constructed as:
%The outputs of the two approximating neural networks for $f[x(k)]$ and $g[x(k)]$ are, respectively, described by

\begin{equation}
	\begin{aligned}
		\hat{f}=\hat{f}[\hat{w}_f, x(k)] = \hat{w}_f^T h_f[x(k)] ,
	\end{aligned}
\end{equation}
\begin{equation}
	\begin{aligned}
		\hat{g}=\hat{g}[\hat{w}_g, x(k)] = \hat{w}_g^T h_g[x(k)] ,
	\end{aligned}
\end{equation}
where $\hat{w}_f$ and $\hat{w}_g$ are vectors representing the neural network output layer parameters, $h_f[x(k)]$ and $h_g[x(k)]$ are Gaussian radial basis function vectors, whose $i$-th element are respectively calculated as
\begin{equation}
	\begin{aligned}
		h_{f_i}[x(k)] = \exp\left\lbrace \frac{-||x(k)-\hat{c}_{f_i}||^2}{2\hat{b}_{f_i}^2}\right\rbrace  ,
	\end{aligned}
\end{equation}
\begin{equation}
	\begin{aligned}
		h_{g_i}[x(k)] = \exp\left\lbrace \frac{-||x(k)-\hat{c}_{g_i}||^2}{2\hat{b}_{g_i}^2}\right\rbrace  ,
	\end{aligned}
\end{equation}
where $\hat{c}_{f_i}$ and $\hat{c}_{g_i}$ are the coordinates of the centre for the $i$-th Gaussian radial basis function, $\hat{b}_{f_i}^2$ and $\hat{b}_{g_i}^2$ are the variances of the $i$-th basis function.

\begin{figure}[htbp]
	\centering
	\includegraphics[width=0.4\textwidth]{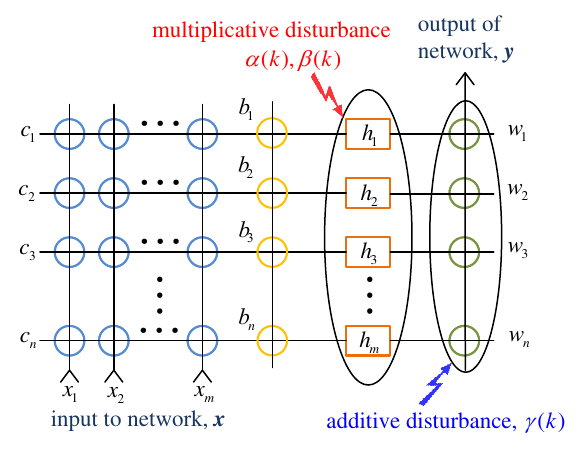}\\
	\caption{The graphic representation of equation~(\ref{equation:specialized neural network}), i.e., the specialized neural network. The blue circle denotes a minus calculation; the yellow circle means a division calculation; the brown rectangle represents the calculation of first minus 2-norm and then exponential; the green circle means a multiplication calculation. $c_{i}=[\hat{c}_{f_i};\hat{c}_{g_i}]$ $b_{i}=[\sqrt{2}\hat{b}_{f_i};\sqrt{2}\hat{b}_{g_i}]$, $h_{i}=[h_{f_{i}}[\cdot];h_{g_{i}}[\cdot] u(k)]$, and $w_{i}=[\hat{w}_{f_{i}}^T; \hat{w}_{g_{i}}^T]$. $\alpha(k)$ and $\beta(k)$ enter the neural network in a multiplicative way, while $\gamma(k)$ enters in an additive way.}\label{Fig_NNdist}
\end{figure}

When taking into consideration the additive and multiplicative disturbances, the unknown nonlinear system affected by disturbances in equation~(\ref{system}) can be approximated by: 
\begin{equation}\label{equation:specialized neural network}
	\begin{aligned}
		y(k+1) \approx & \alpha(k)\hat{w}_f^T h_f[x(k)]+\beta(k)\hat{w}_g^T h_g[x(k)] u(k) \\
		&+\gamma(k)+e(k),\\
	\end{aligned}
\end{equation}
which is actually a neural network structure, as visualized in Fig.~\ref{Fig_NNdist}.
%and Fig.\ref{Fig_NNdist} shows the neural network that corrupted by multiplicative and additive disturbances.

It is a crucial step to recognize the SNN parameters to derive a control law based on the SNN. We hypothesize that the nonlinear dynamics represented by $f[x(k)]$ and $g[x(k)]$ in equation~\ref{system} are stable, while the disturbances represented by $\alpha(k)$ and $\beta(k)$ vary over time. As such, we can identify the SNN parameters in two folds: (i) offline estimate the nonlinear dynamics related parameters i.e., $\hat{c}_{f}$, $\hat{c}_{g}$, $\hat{b}_{f}$, $\hat{b}_{g}$, $\hat{w}_{f}$, and $\hat{w}_{g}$, and (ii) learn the disturbance parameters $\alpha(k)$, $\beta(k)$, and $\gamma(k)$ in real-time.

For the offline parameter estimation, we assume that the multiplicative disturbances $\alpha$ and $\beta$ are equal to a constant, while the additive disturbance $\gamma$ equals to zero. Thus, according to equation~\ref{equation:specialized neural network}, the nonlinear system to recognize turned out to be:
\begin{equation}
	\begin{aligned}
		\hat{y}(k+1) = \hat{w}_f^T h_f[x(k)] + \hat{w}_g^T h_g[x(k)] u(k)
	\end{aligned}
\end{equation}
whose parameters can be realized using a group of historical system measurements. For the online SNN parameter recognition, we need to design a learning approach that allows the estimation of disturbances decoupled from the control law derivation.

\begin{lemma}\label{lemma1}
	Suppose that the approximation value of the multiplicative disturbances and additive disturbance at the $k$-th iteration are known as $\hat{\alpha}(k)$, $\hat{\beta}(k)$, and $\hat{\gamma}(k)$, respectively. According to the dual-property based feedback control law~\cite{fabri1998dual}, the control law can be derived as: 
	\begin{equation}\label{eq:dual control law}
		\begin{aligned}
			u_t(k) =& \frac{[y_r(k+1)-\hat{\alpha}(k)\hat{w}_f h_f-\hat{\gamma}(k)]\hat{\beta}(k)\hat{w}_g h_g}{(1-\lambda)\hat{w}_g h_g P_{\beta} + [\hat{\beta}(k)\hat{w}_g h_g]^2}\\
			&- \frac{(1-\lambda)(\hat{w}_f h_f P_{\alpha\beta}+P_{\gamma\beta})\hat{w}_g h_g}{(1-\lambda)\hat{w}_g h_g P_{\beta} + [\hat{\beta}(k)\hat{w}_g h_g]^2}\\
		\end{aligned}
	\end{equation}
\end{lemma}
by minimizing the sub-optimal dual control performance index
\begin{equation}\label{eq:former index}
	\begin{aligned}
		J(k+1) = & E\{ [y(k+1)-y_r(k+1)]^2 \\
		& - \lambda [y(k+1)-\hat{y}(k+1)]^2|  \mathfrak{I}^k   \}
	\end{aligned}
\end{equation}
instead of the $J(k+1)$ in equation~(\ref{eq_performaceindex}). Note that in equation~(\ref{eq:dual control law}), $\hat{c}_{f_i}$, $\hat{b}_{f_i}$, $\hat{c}_{g_i}$, $\hat{b}_{g_i}$ are hyperparameters in the basis functions $h_f[x(k)$ and $ h_g[x(k)]$, and $\hat{w}_f$, $\hat{w}_g$, which have been trained offline. $P_{\alpha\beta}$, $P_{\gamma\beta}$, and $P_{\beta}$ are the elements of the disturbance estimation error variances, and $\lambda$ is a dual control coefficient. In the following sections, we explain how we estimate the disturbance parameters to achieve the control law derivation.

%\subsection{The finite approximate for multiplicative and additive disturbances}\label{section3.2}
\subsection{The approximation of SNN disturbances via a finite set of candidate values}\label{section3.2}

We design an approach to approximate the additive and multiplicative disturbances as a basis for learning disturbances in the SNN. Our approach aims to avoid identifying the disturbance parameters directly from an unknown disturbance space. Instead, we assign possible candidate values for the disturbance with unknown probability and turn the disturbance learning into probability estimation for those disturbance candidates.

%We design the disturbance learning in SNN to be decoupled from the control law derivation. Our learning approach avoids identifying the disturbance parameters directly from a unknown disturbance space. Instead, we assign possible candidate values for the disturbance with unknown probability, and turn the disturbance estimation into probability learning for those disturbance candidates.

According to assumption~\ref{assumption1}, the disturbances $\alpha(k)$, $\beta(k)$, and $\gamma(k)$ are unknown, time-varying, and yet bounded. Therefore, we can divide those disturbances bounded in the intervals $[\alpha_{l}, \alpha_{u}]$, $[\beta_{l}, \beta_{u}]$, and $[\gamma_{l}, \gamma_{u}]$ into several smaller sub-intervals, and check which sub-interval the disturbance belongs to at a given instant. When the target sub-interval is recognized, we can utilize the middle point of the sub-interval to approximate the true value of the disturbance. Theorem~\ref{thm1} details the partition of sub-intervals for the disturbance intervals.

\begin{thm}\label{thm1}
	Given arbitrary constants $\varepsilon_{\alpha}>0$, $\varepsilon_{\beta}>0$, and $\varepsilon_{\gamma}>0$, there exist positive integers $s_{\alpha}$, $s_{\beta}$, $s_{\gamma}$, and three groups of division points for the disturbance $\alpha$, $\beta$, and $\gamma$:
	\begin{equation}
		\begin{aligned}
			&\alpha_l=\alpha_1<\alpha_2<\cdots<\alpha_{s_{\alpha}}=\alpha_u, \\
			&\beta_l=\beta_1<\beta_2<\cdots<\beta_{s_{\beta}}=\beta_u, \\
			&\gamma_l=\gamma_1<\gamma_2<\cdots<\gamma_{s_{\gamma}}=\gamma_u,
		\end{aligned}
	\end{equation}
	such that
	\begin{equation} \label{thr1.2}
		\begin{aligned}
			&\bigcup_{i=1}^{s_{\alpha}-1}[\alpha_{i},\alpha_{i+1}]=[\alpha_l,\alpha_u], \quad |\alpha_{i}-\alpha_{i+1}|<\varepsilon_{\alpha}, \\
			&\bigcup_{i=1}^{s_{\beta}-1}[\beta_{i},\beta_{i+1}]=[\beta_l,\beta_u], \quad |\beta_{i}-\beta_{i+1}|<\varepsilon_{\beta}, \\
			&\bigcup_{i=1}^{s_{\gamma}-1}[\gamma_{i},\alpha_{i+1}]=[\gamma_l,\gamma_u], \quad |\gamma_{i}-\gamma_{i+1}|<\varepsilon_{\gamma}. \\
		\end{aligned}
	\end{equation}
\end{thm}

%\begin{proof}
\textit{Proof}: For any positive number $\varepsilon_{\alpha}$, there exists a positive integer $s_{\alpha}$ that satisfies
	\begin{equation}\label{eq_p1}
		\begin{aligned}
			s_{\alpha} -1 = \left[ \frac{\alpha_u-\alpha_l}{\varepsilon_{\alpha}} \right] ,
		\end{aligned}
	\end{equation}	
	where $[x]$ represents the maximum integer less than $x$. Then
	\begin{equation}\label{eq_p2}
		\begin{aligned}
			\frac{\alpha_u-\alpha_l}{\varepsilon_{\alpha}} - \left[ \frac{\alpha_u-\alpha_l}{\varepsilon_{\alpha}} \right] <1.
		\end{aligned}
	\end{equation}	
	Submitting equation (\ref{eq_p1}) into equation (\ref{eq_p2}) yields to
	\begin{equation}\label{eq_p3}
		\begin{aligned}
			\frac{\alpha_u-\alpha_l}{\varepsilon_{\alpha}} < \left[ \frac{\alpha_u-\alpha_l}{\varepsilon_{\alpha}} \right] +1=s_{\alpha}.
		\end{aligned}
	\end{equation}	
	The equation (\ref{eq_p3}) is equal to
	\begin{equation}\label{eq_p4}
		\begin{aligned}
			\frac{\alpha_u-\alpha_l}{s_{\alpha}} <  \varepsilon_{\alpha}.
		\end{aligned}
	\end{equation}	
	We divide the interval into $s_{\alpha}$ sub-intervals of the same length. The length of each sub-interval is set as
	\begin{equation}\label{eq_p5}
		\begin{aligned}
			\frac{\alpha_u-\alpha_l}{s_{\alpha}} ,
		\end{aligned}
	\end{equation}	
	then the division points for the interval $[\alpha_{l}, \alpha_{u}]$ are
	\begin{equation}\label{eq_p6}
		\begin{aligned}
			&\alpha_1 = \alpha_l,\\
			&\alpha_i = \alpha_1+\frac{\alpha_u-\alpha_l}{s_{\alpha}}i, \quad i=1,2,...,s_{\alpha}-1\\
			&\alpha_{s_{\alpha}+1} = \alpha_u.\\
		\end{aligned}
	\end{equation}	
	The proof of the division points for the interval $[\beta_{l}, \beta_{u}]$ and $[\gamma_{l}, \gamma_{u}]$ are the same as the process for $[\alpha_l,\alpha_u]$ shown above.
%\end{proof}

Theorem~$\ref{thm1}$ shows that the disturbances $\alpha(k)$, $\beta(k)$, and $\gamma(k)$ can be located in one of the sub-intervals within the disturbance scope. We use the middle point in the sub-interval
\begin{equation}\label{}
	\begin{aligned}
		&\theta_{\alpha_i} = \frac{\alpha_{i-1}-\alpha_i}{2}\\
		&\theta_{\beta_i} = \frac{\beta_{i-1}-\beta_i}{2}\\
		&\theta_{\gamma_i} = \frac{\gamma_{i-1}-\gamma_i}{2}\\
	\end{aligned}
\end{equation}	
to approximate the true disturbances belonging to the sub-intervals $[\alpha_{i-1},\alpha_i]$, $[\beta_{i-1},\beta_i]$, and $[\gamma_{i-1},\gamma_i]$, respectively. Thus, we can use a series of discrete middle points to approximate the whole bounded intervals of the disturbances:
\begin{equation}\label{equation:approximation}
	\begin{aligned}
		&[\alpha_l,\alpha_u] = \Omega_{\alpha} = \{ \theta_{\alpha_1}, \theta_{\alpha_2}, \cdots, \theta_{\alpha_{s_{\alpha}}} \}, \\
		&[\beta_l,\beta_u] = \Omega_{\beta} = \{ \theta_{\beta_1}, \theta_{\beta_2}, \cdots, \theta_{\beta_{s_{\beta}}} \}, \\
		&[\gamma_l,\gamma_u] = \Omega_{\gamma} = \{ \theta_{\gamma_1}, \theta_{\gamma_2}, \cdots, \theta_{\gamma_{s_{\gamma}}} \}. \\
	\end{aligned}
\end{equation}	
and the disturbance identification can be solved as the estimation of the probability of the discrete values in those finite sets. That is, we can designate a prior probability $\pi_{\alpha_i}(1)=1/s_{\alpha}$, $\pi_{\beta_i}(1)=1/s_{\beta}$, and $\pi_{\gamma_i}(1)=1/s_{\gamma}$ as the initial probability for $\theta_{\alpha_i}$, $\theta_{\beta_i}$, and $\gamma_{\alpha_i}$, respectively. The probabilities can be updated over time until converged.

Note that the accuracy of the disturbance approximation grows with the number of the sub-intervals, since more sub-intervals lead to a higher resolution of the probability distribution. However, more sub-intervals increase the computation load of probability estimation for the finite sets. Therefore, we introduce a method to choose an appropriate number of the sub-intervals for users to adjust the trade-off between the disturbance approximation accuracy and the computation, as shown below. 

Suppose that $\alpha^*$ is the true value of the multiplicative disturbance $\alpha(k)$, and it is located in the interval $[\alpha_{i-1},\alpha_i]$. According to equation (\ref{thr1.2}), we can obtain the inequality $|\alpha^*-\theta_{\alpha_i}|<\varepsilon_{\alpha}$ meaning the approximation error is smaller than $\varepsilon_{\alpha}$. We use $\varepsilon_{\alpha}$ as the indicator for adjusting the sub-interval length. I.e., according to equation~(\ref{eq_p1}), we can set the number of the sub-intervals as $s_{\alpha}= \left[ \frac{\alpha_u-\alpha_l}{\varepsilon_{\alpha}} \right] +1$ for the disturbance $\alpha(k)$. Such a sub-interval length aims at meeting a disturbance approximation error specified by users. The choice of sub-interval number for disturbances $\beta(k)$ and $\gamma(k)$ is the same as the case of $\alpha(k)$.

With the disturbance approximation for the multiplicative and additive disturbance, the control objective proposed in Section~$\ref{section2}$ can be realized trough solving the following control problem (P):
\begin{equation}\label{eq:problem_P}
	\begin{aligned}
		&(P)\quad \min J = E \{ [y(k+1)-y_r(k+1)]^2  \\
		& \quad \quad \quad \quad \quad \quad - \lambda [y(k+1)-\hat{y}(k+1)]^2|  \mathfrak{I}^k  \} , \\
		&s.t. \quad y(k+1) = \alpha(k)f[x(k)]+\beta(k)g[x(k)]u(k)\\
		&\quad \quad \quad \quad \quad \quad \quad +\gamma(k)+e(k),\\
	\end{aligned}
\end{equation}	
where the unknown disturbances are within finite sets that: $\alpha(k) \in  \Omega_{\alpha} = \{ \theta_{\alpha_1}, \theta_{\alpha_2}, \cdots, \theta_{\alpha_{s_{\alpha}}} \} $, $\beta(k) \in \Omega_{\beta} = \{ \theta_{\beta_1}, \theta_{\beta_2}, \cdots, \theta_{\beta_{s_{\beta}}} \}$, $\gamma(k) \in \Omega_{\gamma} = \{ \theta_{\gamma_1}, \theta_{\gamma_2}, \cdots, \theta_{\gamma_{s_{\gamma}}} \}$. In the following section, we explain our active learning approach for the disturbance recognition based on the disturbance approximation, and the entailed control law derivation where the disturbance learning and the output controlling are decoupled.

%We designate the prior probability $\pi_{\alpha_i}(1)=1/s_{\alpha}$ as the initial value for the disturbance $\alpha(i) $, which is also the prior probability for $\theta_{\alpha_i}$ in the finite set in equation~(\ref{equation:approximation}). 
%and also it represents that the probability of each $\theta_{\alpha_i}$ that regarded as the true disturbance is equal.

\subsection{Design of the anti-disturbance dual control} \label{section3.3}

We elaborate our anti-disturbance control solution based on the designed SNN and the disturbance approximation. We illustrate our solution from a high-level view and explain the active learning for the disturbances during the control. We present the control law derivation and the active learning together since the two procedures occur simultaneously and iteratively. We explain why our solution decouples the disturbance learning and output tracking control and highlight our solution's properties by proof. We further clarify our anti-disturbance control strategy regarding the efficiency of the active learning and the optimality of the control. 

%This section details the derivation of our proposed anti-disturbance dual control law. 

We derive our control law as follows. Recall that in the control problem $(P)$ in equation~(\ref{eq:problem_P}), the disturbances $\alpha$, $\beta$, and $\gamma$ at the $k$-th iteration is one of the elements in the candidate sets $ \Omega_{\alpha} = \{ \theta_{\alpha_i} : i=1,2, \cdots ,s_{\alpha}\}$, $ \Omega_{\beta} = \{ \theta_{\beta_j} : j=1,2, \cdots ,s_{\beta}\}$, and $ \Omega_{\gamma} = \{ \theta_{\gamma_l} : l=1,2, \cdots ,s_{\gamma}\}$, respectively. For the sake of brevity, we define the disturbance vector as $\theta_t = [\theta_{\alpha_i}, \theta_{\beta_j}, \theta_{\gamma_l}]$, and the disturbance vector set $\Omega$ is
\begin{equation}\label{}
	\begin{aligned}
		\Omega = \left\lbrace  \theta_t : t = 1, 2, \cdots, s_{\alpha}s_{\beta}s_{\gamma}  \right\rbrace ,
	\end{aligned}
\end{equation}
meaning the disturbance vector $\theta_t$ has $s_{\alpha}s_{\beta}s_{\gamma}$ possible options. We define the prior probability for the disturbance vector $\theta_t$ as $\pi(\theta_t)$, and the sum of the prior probabilities is
\begin{equation}\label{}
	\begin{aligned}
		\sum_{t=1}^{s_{\alpha}s_{\beta}s_{\gamma}}\pi(\theta_t) = 1.
	\end{aligned}
\end{equation}
We can express our anti-disturbance dual control law as:
%Our proposed anti-disturbance dual control law is given by 
\begin{equation}\label{controllaw}
	\begin{aligned}
		u(k) = \sum_{t=1}^{s_{\alpha}s_{\beta}s_{\gamma}}\pi(\theta_t|\mathfrak{I}^k) u(k,\theta_t) ,
	\end{aligned}
\end{equation}	
where $\pi(\theta_t|\mathfrak{I}^k)$ is the Bayesian posterior probability for the disturbance vector $\theta_t$ given the information state $\mathfrak{I}^k$. Note that $\pi(\theta_t|\mathfrak{I}^k)$ is updated over time by:% the update recursive equation for this posterior probability is 
\begin{equation} \label{thmposterior}
	\begin{aligned}
		\pi(\theta_t|\mathfrak{I}^k) = \frac{p(y(k)|\theta_t,\mathfrak{I}^k) \pi(\theta_t|\mathfrak{I}^{k-1})}{\sum_{t=1}^{s_{\alpha}s_{\beta}s_{\gamma}}p(y(k)|\theta_t,\mathfrak{I}^k) \pi(\theta_t|\mathfrak{I}^{k-1})} ,
	\end{aligned}
\end{equation}	
where
\begin{equation}\label{thmpdf}
	\begin{aligned}
		p(y(k)|\theta_t,\mathfrak{I}^k) = &\frac{1}{\sqrt{2\pi\Sigma_y(k,\theta_t)}}
		&\exp \left\lbrace -\frac{\tilde{y}^2(k,\theta_t)}{2\Sigma_y(k,\theta_t)} \right\rbrace
	\end{aligned}
\end{equation}	
with the initial posterior probability
\begin{equation}
	\begin{aligned}
		&\pi(\theta_t|\mathfrak{I}^0)=\pi(\theta_t) = \frac{1}{s_{\alpha}s_{\beta}s_{\gamma}}, \\
		&t=1,2, \cdots, s_{\alpha}s_{\beta}s_{\gamma} ,
	\end{aligned}
\end{equation}		
and
\begin{equation}
	\begin{aligned}
		\sum_{t=1}^{ s_{\alpha}s_{\beta}s_{\gamma}}\pi(\theta_t|\mathfrak{I}^k)=1 .
	\end{aligned}
\end{equation}	
$u(k,\theta_t)$ in our control law $(\ref{controllaw})$ represents the feedback control law under a specific combination of disturbance candidates denoted by the disturbance vector $\theta_t$. According to Lemma~\ref{lemma1}, the feedback control law corresponding to the disturbance $\theta_t$, which we call the $\theta_t$-candidate feedback control law, can be written as:
\begin{equation}\label{nconditionuk}
	\begin{aligned}
		u_t(k) =& \frac{[y_r(k+1)-\theta_t(1)\hat{w}_f h_f-\theta_t(3)]\hat{\beta}(k)\hat{w}_g h_g}{(1-\lambda)\hat{w}_g h_g P_{\beta} + [\theta_t(2)\hat{w}_g h_g]^2}\\
		&- \frac{(1-\lambda)(\hat{w}_f h_f P_{\alpha\beta}+P_{\gamma\beta})\hat{w}_g h_g}{(1-\lambda)\hat{w}_g h_g P_{\beta} + [\theta_t(2)\hat{w}_g h_g]^2}\\
	\end{aligned}
\end{equation}
where the $\theta_t(1)$, $\theta_t(2)$, and $\theta_t(3)$ is the first, second, and third element of the $\theta_t$, meaning the approximation of disturbance $\alpha(k)$, $\beta(k)$, and $\gamma(k)$, respectively. $P_{\beta}$, $P_{\alpha\beta}$, and $P_{\gamma\beta}$ are the elements of the disturbance vector estimation error covariance matrix $P(k)=E\{\tilde{\theta}_t\tilde{\theta}_t^T\}$. We design the update of the disturbance vector estimation error covariance matrix as:
\begin{equation}\label{varianceP}
	\begin{aligned}
		P(k+1) = \log_2 [\eta/\pi(\theta_t|\mathfrak{I}^k)+1] P(k),
	\end{aligned}
\end{equation}	
where the coefficient $\eta$ is set as $1/s_{\alpha}s_{\beta}s_{\gamma}$. The design of this update aims to assign higher covariance for the element with lower posterior probability.

\begin{figure}[htbp]
	\centering
	\includegraphics[width=0.5\textwidth]{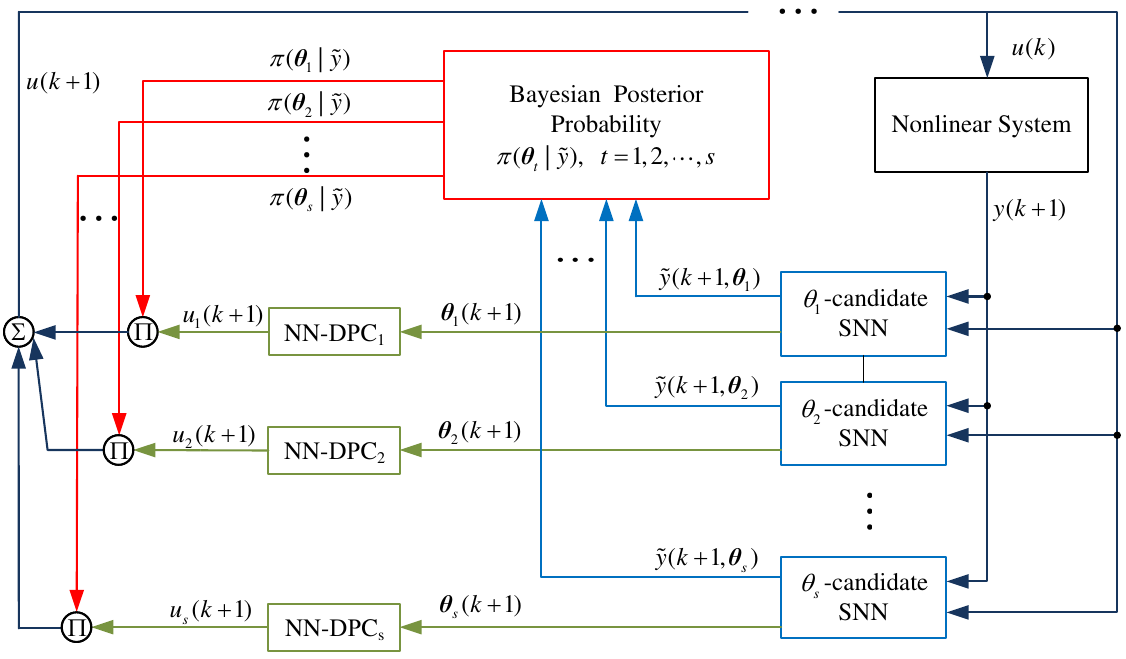}\\
	\caption{The block diagram of the anti-disturbance dual control scheme. The symbol ``$\Pi$'' represents the multiplication operator and the symbol ``$\Sigma$'' is the adding operator.}\label{Fig_1}
\end{figure}

Fig.\ref{Fig_1} visualizes our control strategy by a block diagram. Following equation~(\ref{controllaw}), the designed control law is constructed by summing up $s(s=s_{\alpha}s_{\beta}s_{\gamma})$ $\theta_t$-candidate control laws weighted by the corresponding Bayesian posterior probabilities $\pi(\theta_t|\tilde{y})$, as shown in the diagram. This control strategy holds the conflicting dual properties: (i) driving the system output to track the reference trajectory and (ii) performing active learning of system dynamics and disturbances estimation to reduce the uncertainties for the control. As shown in the diagram, we decouple the tracking control and the active learning during the control process. Specifically, we use a special neural network specialized by the disturbance $\theta_t\ (t\in \{1,2,3,...,s_{\alpha}s_{\beta}s_{\gamma}\})$, which we call the $\theta_i$-candidate SNN, to derive the corresponding dual-property control law (NN-DPC$_t$). We achieve the active learning of the unknown disturbances $\theta_t$ by updating the Bayesian posterior probability of the disturbance candidates, which is calculated based on the one-step-ahead output prediction error $\tilde{y}(k+1,\theta_t)$ in the $\theta_t$-candidate SNN. Hence the derivation of tracking control law for each candidate SNN and the active learning for the Bayesian posterior probability of the corresponding disturbance candidate are separated. Thus, the disturbance learning is decoupled from the influence of control law derivation.

In that way, our anti-disturbance dual control holds the following feature. Suppose the true value of the disturbances is mostly close to the candidate disturbance $\theta_{t^{*}}$, then the learned disturbance can converge to $\theta_{t^{*}}$ by our active learning, and the entailed $\theta_{t^{*}}$-candidate one-step-ahead prediction $\hat{y}(k+1,\theta_{t^{*}})$ is the optimal system output estimation under our approach. We clarify and provide proof for this statement in the following two theorems. Particularly, Theorem~\ref{thminequality} proves that the variance of the $\theta_{t^{*}}$-candidate one-step-ahead prediction error  $\tilde{y}(k+1,\theta_{t^{*}})$ is the least compared with $\tilde{y}(k+1,\theta_t), t \neq t^{*}$. Theorem~\ref{thm2} provides the proof for the convergence of the disturbance candidate's Bayesian posterior probability. I.e., the posterior probability of the $\theta_{t^{*}}$-candidate $\pi(\theta_{t^{*}}|\tilde{y}) \rightarrow 1$, when the $k\rightarrow \infty$, and the other posterior probability $\pi(\theta_t|\tilde{y}) \rightarrow 0$ for $\forall t \neq t^{*} $. 

\begin{thm}\label{thminequality}
	Suppose the disturbance vector $\theta_{t^{*}}$ is most close to the true value for $\alpha(k)$, $\beta(k)$, and $\gamma(k)$. Then for any disturbance vectors $ \theta_t, t \neq t^{*}$, we have the inequality:
	\begin{equation}\label{sigma_inq}
		\begin{aligned}
			E\left\lbrace \tilde{y}^2(k+1|k,\theta_t) \right\rbrace > E\left\lbrace \tilde{y}^2(k+1|k,\theta_{t^{*}}) \right\rbrace.
		\end{aligned}
	\end{equation}
	where $E\left\lbrace \tilde{y}^2(k+1|k,\theta_t) \right\rbrace$ represents the variance of error of the system output estimation when using $\theta_t$ for the estimation at the instant $k$.
\end{thm}

%\begin{proof}
\textit{Proof}:	We define the one-step-ahead prediction of the system output at the $k$-th iteration using the disturbance vector $\theta_{t^{*}}$ as 
	\begin{equation}\label{}
		\begin{aligned}
			\hat{y}(k+1,\theta_{t^{*}}) = &\theta_{t^{*}}(1) \hat{f}[x(k)] + \theta_{t^{*}}(2) \hat{g}[x(k)] u(k)\\
			& + \theta_{t^{*}}(3) ,
		\end{aligned}
	\end{equation}
	where $\theta_{t^{*}}(1)$, $\theta_{t^{*}}(2)$, and $\theta_{t^{*}}(3)$ are the first, second, and third element of the $\theta_{t^{*}}$, respectively. The one-step-ahead prediction error is defined as
	\begin{equation}\label{}
		\begin{aligned}
			\tilde{y}(k+1,\theta_{t^{*}}) = y(k+1) - \hat{y}(k+1,\theta_{t^{*}})
		\end{aligned}
	\end{equation}
	Submitting the system described in equation~(\ref{system}) into the $\tilde{y}(k+1,\theta_{t^{*}})$ yields
	\begin{equation}\label{}
		\begin{aligned}
			\tilde{y}(k+1,\theta_{t^{*}}) = &\{ \alpha(k)f[x(k)]-\theta_{t^{*}}(1) \hat{f}[x(k)] \} \\
			&+\left\lbrace \beta(k)g[x(k)]-\theta_{t^{*}}(2) \hat{g}[x(k)]\right\rbrace  u(k) \\
			&+ \{\gamma(k)-\theta_{t^{*}}(3) \} +e(k)
		\end{aligned}
	\end{equation}
	Suppose the learned nonlinear functions are very close to the true ones and satisfy $f[x(k)]=\hat{f}[x(k)]$, $g[x(k)]=\hat{g}[x(k)]$. Then, the one-step-ahead prediction error under $\theta_{t^*}$ can be rewritten as 
	\begin{equation}\label{predictionerror2}
		\begin{aligned}
			\tilde{y}(k+1, \theta_{t^*}) & = \tilde{\alpha}_{t^*}(k)\hat{f}[x(k)]+\tilde{\beta}_{t^*}(k)\hat{g}[x(k)]u(k)\\
			&\quad +\tilde{\gamma}_{t^*}(k)+e(k)\\
			& = \tilde{\theta}_{t^*}(k) \phi(k) + e(k).
		\end{aligned}
	\end{equation}
	where $\tilde{\alpha}_{t^*}(k)=\alpha(k)-\theta_{t^{*}}(1)$, $\tilde{\beta}_{t^*}(k)=\beta(k)-\theta_{t^{*}}(2)$, $\tilde{\gamma}_{t^*}(k)-\theta_{t^{*}}(3)$, $\phi(k)=[\hat{f}[x(k)], \hat{g}[x(k)]u(k), 1]$. We can define the variance of one-step-ahead prediction of the system output under $\theta_{t^{*}}$ as:
	\begin{equation}\label{}
		\begin{aligned}
			&\Sigma_y(k+1,\theta_{t^{*}}) = E\left\lbrace \tilde{y}^2(k+1,\theta_{t^{*}}) \right\rbrace \\
			&=  \phi(k)^T E \left\lbrace \tilde{\theta}_{t^*}(k) \tilde{\theta}_{t^*}^T(k) \right\rbrace \phi(k) +E \left\lbrace e^2(k) \right\rbrace \\
			&=\phi(k)^T P(k) \phi(k) + \sigma^2 .
		\end{aligned}
	\end{equation}
	From assumption \ref{assumption3}, the process error $e(k)$ is a Gaussian process with $0$ mean and $\sigma^2$ variance. Since $\theta_{t^*}$ is most close to the ground truth disturbance and the Bayesian posterior probability for the $\theta_{t^*}$-candidate $\pi(\theta_{t^{*}}|\mathfrak{I}^k)$ converges to 1, the part $\log_2 [\eta/\pi(\theta_t|\mathfrak{I}^k)+1]$ in equation~(\ref{varianceP}) is less than $1$, which leads the estimation error covariance $P(k)$ converges to $0$. Then we get
	\begin{equation}\label{sigma_t}
		\begin{aligned}
			\Sigma_y(k+1,\theta_{t^{*}}) = \sigma^2 .
		\end{aligned}
	\end{equation}
	
	If the disturbance vector $\theta_t$ is not the most close to true value of the disturbance (or $t\neq t^{*}$), then the Bayesian posterior probability $\pi(\theta_t|\mathfrak{I}^k)$ converges to $0$, and the part $\log_2 [\eta/\pi(\theta_t|\mathfrak{I}^k)+1]$ in equation~(\ref{varianceP}) is larger than $1$, which makes the estimation error covariance $P(k)$ much larger than $0$. Therefore, we can obtain the variance of one-step-ahead prediction of the system output under $\theta_t\ (t\neq t_{*})$ as:
	\begin{equation}\label{sigma_i}
		\begin{aligned}
			&\Sigma_y(k+1,\theta_t) = E\left\lbrace \tilde{y}^2(k+1,\theta_t) \right\rbrace\\
			&	= \phi(k)^T P(k) \phi(k) + \sigma^2 > \sigma^2
		\end{aligned}
	\end{equation}
	Compare equation (\ref{sigma_t}) with (\ref{sigma_i}), and we can obtain the inequality that
	\begin{equation}\label{sigma_inq}
		\begin{aligned}
			E\left\lbrace \tilde{y}^2(k+1|k,\theta_t) \right\rbrace > E\left\lbrace \tilde{y}^2(k+1|k,\theta_{t^{*}}) \right\rbrace.
		\end{aligned}
	\end{equation}
	
%\end{proof}

\begin{thm}\label{thm2}
	Suppose $ \theta_{t^{*}}$ is the true value of the disturbance vector to be identified. Then with our active learning approach, we have that the posterior probability $\pi(\theta_{t^{*}}|\mathfrak{I}^k) \rightarrow 1$ when the $k\rightarrow \infty$, while the posterior probability $\pi(\theta_t|\mathfrak{I}^k) \rightarrow 0$ when $k\rightarrow \infty$ for $\forall t \neq t^{*} $.	
\end{thm}

%\begin{proof}
\textit{Proof}:	Supposed the sequence $	\tilde{y}(k|k,\theta_t)$ is an asymptotic weak stationary sequence. The stationary process $\tilde{y}(k|k,\theta_t)$ is ergodic for each $\theta_t$, and the variance $\Sigma_y(k|k,\theta_t)$ is constant that it can be rewritten as $\Sigma_t$.
	
	We define
	\begin{equation} \label{proof_L_def}
		\begin{aligned}
			L_t(k) = \pi(\theta_t|\mathfrak{I}^k) \pi^{-1}(\theta_{t^{*}}|\mathfrak{I}^k) ,\\
			t \in \{1,2,\cdots, s_{\alpha}s_{\beta}s_{\gamma}\}.
		\end{aligned}
	\end{equation}
	Substituting equation (\ref{thmposterior}) into the expression of $L_i(k)$ yields
	\begin{equation} \label{proof_Li}
		\begin{aligned}
			L_t(k) &= \frac{ p(y(k)|\theta_t,\mathfrak{I}^k) \pi(\theta_t|\mathfrak{I}^{k-1}) }{ p(y(k)|\theta_{t^{*}},\mathfrak{I}^k) \pi(\theta_{t^{*}}|\mathfrak{I}^{k-1}) } \\
			& = \frac{p(y(k)|\theta_t,\mathfrak{I}^k)}{p(y(k)|\theta_{t^{*}},\mathfrak{I}^k)} L_t(k-1)\\
		\end{aligned}
	\end{equation}
	According to equation (\ref{proof_Li}), we can obtain that
	\begin{equation} \label{proof_Li_n}
		\begin{aligned}
			\frac{L_t(k+n-1) }{L_t(k-1) } = \prod_{\tau=k}^{k+n-1} \frac{p(y(\tau)|\theta_t,\mathfrak{I}^{\tau})}{p(y(\tau)|\theta_{t^{*}},\mathfrak{I}^{\tau})},
		\end{aligned}
	\end{equation}
	Submitting the probability density function (\ref{thmpdf}) into equation (\ref{proof_Li_n}) yields
	\begin{equation}
		\begin{aligned}
			&\frac{L_t(k+n-1) }{L_t(k-1) }
			= \left\lbrace \frac{\Sigma_t}{\Sigma_{t^{*}}} \right\rbrace ^{-\frac{n}{2}}
			\frac{ \exp\left\lbrace \sum_{\tau=k}^{k+n-1} -\frac{\tilde{y}^2(\tau,\theta_t)}{2\Sigma_t} \right\rbrace }  {\exp\left\lbrace \sum_{\tau=k}^{k+n-1} -\frac{\tilde{y}^2(\tau,\theta_{t^{*}})}{2\Sigma_{t^{*}}}\right\rbrace }
		\end{aligned}
	\end{equation}
	The natural logarithm of $\frac{L_t(k+n-1) }{L_t(k-1) } $ can be obtained as:
	\begin{equation} \label{proof_ln_L}
		\begin{aligned}
			&\ln\left\lbrace \frac{L_t(k+n-1) }{L_t(k-1) }  \right\rbrace = \frac{n}{2} \ln \left\lbrace \frac{\Sigma_{t^{*}}}{\Sigma_t} \right\rbrace \\
			&\quad - \frac{1}{2} \sum_{\tau=k}^{k+n-1} \frac{\tilde{y}^2(\tau|\tau,\theta_t)}{\Sigma_t}
			+ \frac{1}{2} \sum_{\tau=k}^{k+n-1}\frac{\tilde{y}^2(\tau|\tau,\theta_{t^{*}})}{\Sigma_{t^{*}}}
		\end{aligned}
	\end{equation}
	With the assumption that $ \theta_{t^{*}}$ is the true value of the disturbance vector, and the process $\tilde{y}(k,\theta_{t^{*}})$ is ergodic, we can obtain
	\begin{equation}
		\begin{aligned}
			E\left\lbrace \tilde{y}^2(k|k,\theta_{t^{*}}) \right\rbrace = \lim_{n \rightarrow \infty } \frac{1}{n} \sum_{\tau=k}^{k+n-1} \tilde{y}^2(\tau|\tau,\theta_{t^{*}})= \Sigma_{t^{*}}
		\end{aligned}
	\end{equation}
	the right of which can be also written as:
	\begin{equation} \label{proof_lim_1}
		\begin{aligned}
			\lim_{n \rightarrow \infty } \frac{1}{n} \sum_{\tau=k}^{k+n-1}\frac{\tilde{y}^2(\tau|\tau,\theta_{t^{*}})}{\Sigma_{t^{*}}} = 1
		\end{aligned}
	\end{equation}
	While for any other disturbance vector $ \theta_t\ (t \neq t^{*})$ , we get
	\begin{equation}
		\begin{aligned}
			E\left\lbrace \tilde{y}^2(k|k,\theta_t) \right\rbrace = \lim_{n \rightarrow \infty } \frac{1}{n} \sum_{\tau=k}^{k+n-1} \tilde{y}^2(\tau|\tau,\theta_t),
		\end{aligned}
	\end{equation}
	and according to the inequality (\ref{sigma_inq}), we get the limitation that
	\begin{equation} \label{proof_lim_2}
		\begin{aligned}
			\lim_{n \rightarrow \infty } \frac{1}{n} \sum_{\tau=k}^{k+n-1} \tilde{y}^2(\tau|\tau,\theta_t) > \Sigma_{t^{*}} , \quad \forall t \neq t^{*}
		\end{aligned}
	\end{equation}
	Submitting equations (\ref{proof_lim_1}) and (\ref{proof_lim_2}) into (\ref{proof_ln_L}) yields
	\begin{equation} \label{proof_ln_L2}
		\begin{aligned}
			&\lim_{n \rightarrow \infty }\frac{2}{n}\ln\left\lbrace \frac{L_t(k+n-1) }{L_t(k-1) }  \right\rbrace =  \ln \left\lbrace \frac{\Sigma_{t^{*}}}{\Sigma_t} \right\rbrace \\
			&- \lim_{n \rightarrow \infty } \frac{1}{n} \sum_{\tau=k}^{k+n-1} \frac{\tilde{y}^2(\tau|\tau,\theta_t)}{\Sigma_t}
			+ \lim_{n \rightarrow \infty } \frac{1}{n} \sum_{\tau=k}^{k+n-1}\frac{\tilde{y}^2(\tau|\tau,\theta_{t^{*}})}{\Sigma_{t^{*}}} \\
			& =  \ln \left\lbrace \frac{\Sigma_{t^{*}}}{\Sigma_t} \right\rbrace - \frac{\Sigma_{t^{*}}}{\Sigma_t}- \frac{M_t}{\Sigma_t} +1,
		\end{aligned}
	\end{equation}
	where
	\begin{equation} \label{}
		\begin{aligned}
			M_t = \lim_{n \rightarrow \infty } \frac{1}{n} \sum_{\tau=k}^{k+n-1} \tilde{y}^2(\tau|\tau,\theta_t)-\Sigma_{t^{*}}.
		\end{aligned}
	\end{equation}
	It is obviously that for $\forall x>0$, we have $ln(x)-x+1\leq0$. Therefore, since in equation~\ref{proof_ln_L2} there is $\Sigma_{t^{*}}>0$ and $\Sigma_t>0$, we can obtain the inequality that
	\begin{equation} \label{proof_ln_ieq1}
		\begin{aligned}
			\ln \left\lbrace \frac{\Sigma_{t^{*}}}{\Sigma_t} \right\rbrace - \frac{\Sigma_{t^{*}}}{\Sigma_t} +1\leq0
		\end{aligned}
	\end{equation}
	According to equation (\ref{proof_lim_2}), we get $M_t>0$. Submitting inequality (\ref{proof_ln_ieq1}) and $M_t>0$ into the limitation (\ref{proof_ln_L2}) yields
	\begin{equation} \label{proof_ln_L3}
		\begin{aligned}
			&\lim_{n \rightarrow \infty }\frac{2}{n}\ln\left\lbrace \frac{L_t(k+n-1) }{L_t(k-1) }  \right\rbrace =  -c_1
		\end{aligned}
	\end{equation}
	where $c_1>0$. According to equation (\ref{proof_ln_L3}), we obtain
	\begin{equation} \label{proof_ln_L4}
		\begin{aligned}
			\lim_{n \rightarrow \infty } L_t(k+n-1)  & =\lim_{n \rightarrow \infty }  c_2 L_t(k-1) \exp\{-nc_1/2\}\\
			&	=0,   \quad \forall t \neq t^{*} ,
		\end{aligned}
	\end{equation}
	where $c_2$ is a constant. Combing equation~(\ref{proof_L_def}) and the limitation~(\ref{proof_ln_L4}), we obtain
	\begin{equation} \label{proof_ln_L5}
		\begin{aligned}
			\lim_{n \rightarrow \infty } \pi(\theta_t|\mathfrak{I}^k) \pi^{-1}(\theta_{t^{*}}|\mathfrak{I}^k) =0,  \quad \forall t \neq t^{*} ,
		\end{aligned}
	\end{equation}
	which equals to
	\begin{equation} \label{proof_ln_L6}
		\begin{aligned}
			\lim_{n \rightarrow \infty }  \pi(\theta_t|\mathfrak{I}^k)=0, \quad \forall t \neq t^{*} .
		\end{aligned}
	\end{equation}
	Since the sum of all the probability of the possible disturbance vector equals $1$, we can get
	\begin{equation} \label{proof_ln_L7}
		\begin{aligned}
			\lim_{n \rightarrow \infty }  \pi(\theta_{t^{*}}|\mathfrak{I}^k)=1 .
		\end{aligned}
	\end{equation}
	indicating that by our active learning approach, the identified disturbance converges to the true value over time. As such, our designed control law $u(k) = \sum_{t=1}^{s_{\alpha}s_{\beta}s_{\gamma}}\pi(\theta_t|\mathfrak{I}^k) u_t(k)$ takes the right disturbance vector over time, and derives an optimal control strategy for the nonlinear system.
%\end{proof}

\begin{remark}
	The proof for theorem \ref{thm2} shows that if the $\theta_{t^{*}}$ is the best approximation for the true disturbance in the infinite set $\Omega = \left\lbrace  \theta_t : t = 1, 2, \cdots, s_{\alpha}s_{\beta}s_{\gamma}  \right\rbrace$, then the posterior probability $\pi(\theta_{t^{*}}|\mathfrak{I}^k) $ converges exponentially from the initial value $1/ s_{\alpha}s_{\beta}s_{\gamma}$ to $1$. Recall that in the definition of $\theta_t=[\theta_{\alpha_i}, \theta_{\beta_j}, \theta_{\gamma_l}]$, $\theta_t$ consists of all the possible composite of the disturbances intervals. Therefore, the convergence of posterior probability also indicates that the true disturbance $\alpha(k)$, $\beta(k)$, and $\gamma(k)$ is located in the sub-intervals $[\alpha_{i-1},\alpha_{i}]$, $[\beta_{j-1},\beta_{j}]$, and $[\gamma_{l-1},\gamma_{l}]$, respectively.
\end{remark}

\begin{remark}
	Equation (\ref{controllaw}) indicates that our controller consists of $s\ (s=s_{\alpha}s_{\beta}s_{\gamma})$ sub-controllers, each of which corresponds to a specific $\theta_t\ (t=1,2, \cdots, s_{\alpha}s_{\beta}s_{\gamma})$, and the control law is expressed as the weighted average of the $\theta_t$-candidate control $u(k,\theta_t)$ weighted by the posterior probabilities $\pi(\theta_t|\mathfrak{I}^k)$. This strategy makes sense in that the posterior probability of the true disturbance $\theta_{t^{*}}$ will converge to $1$, resulting in the $\theta_{t^{*}}$-candidate control law $u(k,\theta_{t^{*}})$ being optimal for the considered nonlinear system.
\end{remark}

\begin{remark}
	Suppose at an instant $k$, the $\theta_{t^{*}}$-candidate posterior probability $\pi(\theta_{t^{*}}|\mathfrak{I}^k) $ has converged to $1$, and the $\theta_t$-candidate posterior probability $\pi(\theta_t|\mathfrak{I}^k) $ equals $0$ $(t\neq t^{*})$. According to equation (\ref{thmposterior}), the posterior probability $\pi(\theta_{t^{*}}|\mathfrak{I}^k) $ and $\pi(\theta_t|\mathfrak{I}^k)\ (t\neq t^{*}) $ will not change and are locked, even the true disturbances $\alpha(k)$, $\beta(k)$ and $\gamma(k)$ vary a lot after the instant $k$. Therefore, the control law is locked at $u(k) \equiv u(k,\theta_{t^{*}})$ from the instant $k$. This is problematic, and it is necessary to design a criterion to check whether the disturbance is changed that the true value is no longer approximated to $\theta_{t^{*}}$ after the instant $k$. In this case, we set such a criterion as follows. If the system output estimation error $\tilde{y}(k+1,\theta_{t^{*}})>\epsilon$ during the stationary process (where one of the posterior probabilities of $\theta_t \ (t\in\{1,2,3,...,s_{\alpha}s_{\beta}s_{\gamma}\})$ is equal to 1), it means the disturbances has changed, where $\epsilon$ is the system output tracking admissible error. When it is detected that the disturbance has changed, we will reset the posterior probabilities as the initial value $1/ s_{\alpha}s_{\beta}s_{\gamma}$ for all $\theta_t \ (t\in\{1,2,3,...,s_{\alpha}s_{\beta}s_{\gamma}\})$ and restart the active learning procedures.
\end{remark}

\section{Simulations and results} \label{section4}
\begin{table}[H]\label{procedure}
	\centering
	\caption{Implementation of the designed approach.}
	\resizebox{8cm}{!}
	{\begin{tabular}{l}
			\toprule
			\textbf{Initialization:}\\
			Initialize the coefficients $\varepsilon_{\alpha}>0$, $\varepsilon_{\beta}>0$, $\varepsilon_{\gamma}>0$, and $\epsilon>0$, \\
			the disturbance vector finite set $\Omega = \left\lbrace  \theta_t : t = 1, 2, \cdots, s_{\alpha}s_{\beta}s_{\gamma}  \right\rbrace$, \\
			the state information $\mathfrak{I}^0$, the Bayesian posterior probability $\pi(\theta_t|\mathfrak{I}^0)$, \\
			the dual-property coefficient $0<\lambda<1$, the disturbance estimation\\
			error covariance matrix $P(0)$, the control $u^*(1)$, and a posterior\\
			probability threshold $\phi$.\\
			\midrule
			\textbf{Computation:}\\
			(1) Apply the control law $u^*(k)$ for the system and observe the \\
			\quad \quad system output $y(k+1)$;\\
			(2) Calculate the Bayesian posterior probability $\pi(\theta_t|\mathfrak{I}^k)$ based  \\
			\quad \quad on equations (\ref{thmposterior}) and (\ref{thmpdf}) at the iteration $k+1$;  \\
			(3) Calculate the $\theta_t$-candidate feedback control law according to   \\
			\quad \quad equation (\ref{nconditionuk}) at the iteration $k+1$;\\
			(4) Calculate our proposed control law $u^*(k+1)$ according to \\
			\quad \quad equation (\ref{controllaw});\\
			(5) Calculate the prediction of the system output at the $(k+1)$-th\\
			\quad \quad iteration $\hat{y}(k+1,\theta_t)$, and the prediction error $\tilde{y}(k+1,\theta_t)$;\\
			(6) \textbf{If} there is $\tilde{y}(k+1,\theta_t)>\epsilon$ and $\pi(\theta_t|\mathfrak{I}^k)>\phi$\\
			\quad \quad \quad \quad The Bayesian posterior probability is reset as   \\
			\quad \quad \quad \quad $\pi(\theta_t|\mathfrak{I}^k)=1/ s_{\alpha}s_{\beta}s_{\gamma}$.\\
			\quad \; \textbf{End if}\\
			(7) Update the disturbance estimation error covariance matrix \\
			\quad \quad $P(k+1)$according to equation (\ref{varianceP})\\
			(8) Repeat steps (1) to (7) in future iterations.\\
			\midrule
			\textbf{Notice:}\\
			(1) The parameters other than the disturbances in the specialized\\
			neural network are trained offline. The training data is collected\\
			in a stationary control process to guarantee the accuracy of\\
			offline parameter identification. \\
			(2) The upper and lower bounds for the multiplicative and additive  \\
			disturbances can be estimated empirically based on the historical data\\
			collected from the control process.\\
			\bottomrule
	\end{tabular}}
	\label{table1}
\end{table}

This section demonstrates the simulation and evaluation of the designed approach for the control of unknown nonlinear systems with diversities of settings. Section~\ref{section4.1} simulates our approach in controlling a numerical nonlinear system that suffers from multiplicative disturbances varying over time, with the reference trajectory being a continuous sine wave for the output tracking control. Section~\ref{section4.2} shows simulations on the same numerical system yet corrupted by additive disturbances, with the reference trajectory being a step featured square wave. Section~\ref{section4.3} further tests and analyzes our approach on the mentioned numerical system via Monte Carlo statistical simulations. In Section~\ref{section4.4}, we evaluate our approach in the speed control of CRH3 high-speed train in the existence of unknown multiplicative and additive disturbances varying over time. We also compare our approach with the ideal benchmark optimal control and a latest model-free adaptive control method in this evaluation.

The implementation procedures for the proposed method are summarized in Table~\ref{procedure}. We run our experiments using MATLAB 2022a on Windows 11, with 2.80 GHz core 16.0 GB RAM. 

\subsection{Case 1: simulations where the nonlinear system is disturbed by multiplicative disturbances and the reference trajectory is a sine wave} \label{section4.1}
Consider an affine nonlinear system \cite{fabri1998dual} as
\begin{equation}\label{case1}
	\begin{aligned}
		&y(k+1) = \alpha(k)[\sin(x(k))+\cos(3x(k))] \\
		&+\beta(k)(2+\cos(x(t)))u(k)+\gamma(k)+e(k+1)
	\end{aligned}
\end{equation}
where the multiplicative disturbances $\alpha(k)$ and $\beta(k)$ vary in the bounded interval $[0.75,1.25]$ and $[0.75,1.05]$, respectively, and the additive disturbances $\gamma(k)$ is zero. The state $x(k)$ is $x(k)=y(k)$, and the system noise $e(k)$ follows Gaussian distribution $N(0,0.0004)$. This system has two nonlinear functions $f[x(k)]=\sin(x(k))+\cos(3x(k))$ and $g[x(k)]=2+\cos(x(t))$. The reference trajectory for output tracking is set as a cosine wave $y_r(k)=\cos(5\pi k/600)$. The true value of the multiplicative disturbances $\alpha(k)$ and $\beta(k)$ are set over time as
\begin{equation} \label{case1mulalpha}
	\begin{aligned}
		\alpha(k) = \left\{ \begin{array}{rcl}
			1.0 & \mbox{for}& 1 \leq k < 85\\
			1.11 & \mbox{for} & 85 \leq k < 180\\
			0.78 & \mbox{for} & 180 \leq k < 340\\
			0.91 & \mbox{for} & 340 \leq k < 520\\
			1.18 & \mbox{for} & 520 \leq k \leq 600\\
		\end{array}\right.
	\end{aligned}
\end{equation}
\begin{equation} \label{case1mulbeta}
	\begin{aligned}
		\beta(k) = \left\{ \begin{array}{rcl}
			0.9 & \mbox{for}& 1 \leq k < 180\\
			0.82 & \mbox{for} & 180\leq k < 340\\
			1.0 & \mbox{for} & 340\leq k \leq 600\\
		\end{array}\right.
	\end{aligned}
\end{equation}
where $k$ denotes instant from 1 to 600. We firstly construct a specialized neural network (SNN) for this nonlinear system, where we obtain parameters for the two nonlinear functions $f[x(k)]$ and $g[x(k)]$ through offline learning. With historical system measurement, $f[x(k)]$ is estimated as the nonlinear function $\hat{f}[x(k)]$ in the SNN that uses the Gaussian radial basis with the center points $\hat{c}_f = [-2, -1.5, -1, -0.5, 0, 0.5, 1, 15, 2]$, the width $\hat{b}_f^2=1$, and the weights $\hat{w}_f = [11.2856, -4.6174, -12.3754, 1.5622,\\ 12.0864, 2.2351, -11.9197, -4.4981, 12.6531]$. Similarly, $g[x(k)]$ is estimated as $\hat{g}[x(k)]$ that has the center points $\hat{c}_g=[-2,  0,  2]$, the variance $\hat{b}_g^2=3.6$, and the weights $\hat{w}_f=[ -0.4449, 3.5097, -0.4449 ]$.

According to theorem \ref{thm1}, we set the admissible disturbance approximation error as $\varepsilon_{\alpha}=0.1$ and $\varepsilon_{\beta}=0.1$. The corresponding sub-interval number can be obtained as $s_{\alpha}=5$ and $s_{\beta}=3$. The multiplicative disturbance approximation finite sets are $ \Omega_{\alpha} = \{ 0.8, 0.9, 1.0, 1.1, 1.2 \}$ and $\Omega_{\beta} = \{ 0.8, 0.9, 1.0 \}$. The number of elements is 15 in the finite disturbance vector set, and the elements in this set are $\theta_1=[0.8, 0.8, 0]$, $\theta_2=[0.8, 0.9, 0]$, $\theta_3=[0.8, 1.0, 0]$, $\theta_4=[0.9, 0.8, 0]$, $\theta_5=[0.9, 0.9, 0]$, $\theta_6=[0.9, 1.0, 0]$, $\theta_7=[1.0, 0.8, 0]$, $\theta_8=[1.0, 0.9, 0]$, $\theta_9=[1.0, 1.0, 0]$, $\theta_{10}=[1.1, 0.8, 0]$, $\theta_{11}=[1.1, 0.9, 0]$, $\theta_{12}=[1.1, 1.0, 0]$, $\theta_{13}=[1.2, 0.8, 0]$, $\theta_{14}=[1.2, 0.9, 0]$, $\theta_{15}=[1.2, 1.0, 0]$. The initial Bayesian posterior probabilities are set as $\pi(\theta_t|\mathfrak{I}^0)=1/15$ for $t=1,2,\cdots, 15$. The unlock criterion parameter is set as $\epsilon=0.08$. The initial disturbance estimation error covariance is set as $P(1)=I_{3\times3}$, and the dual-property control coefficient is set as $\lambda=0.9$.

\begin{figure}[htbp]
	\centering
	\includegraphics[width=0.8\textwidth]{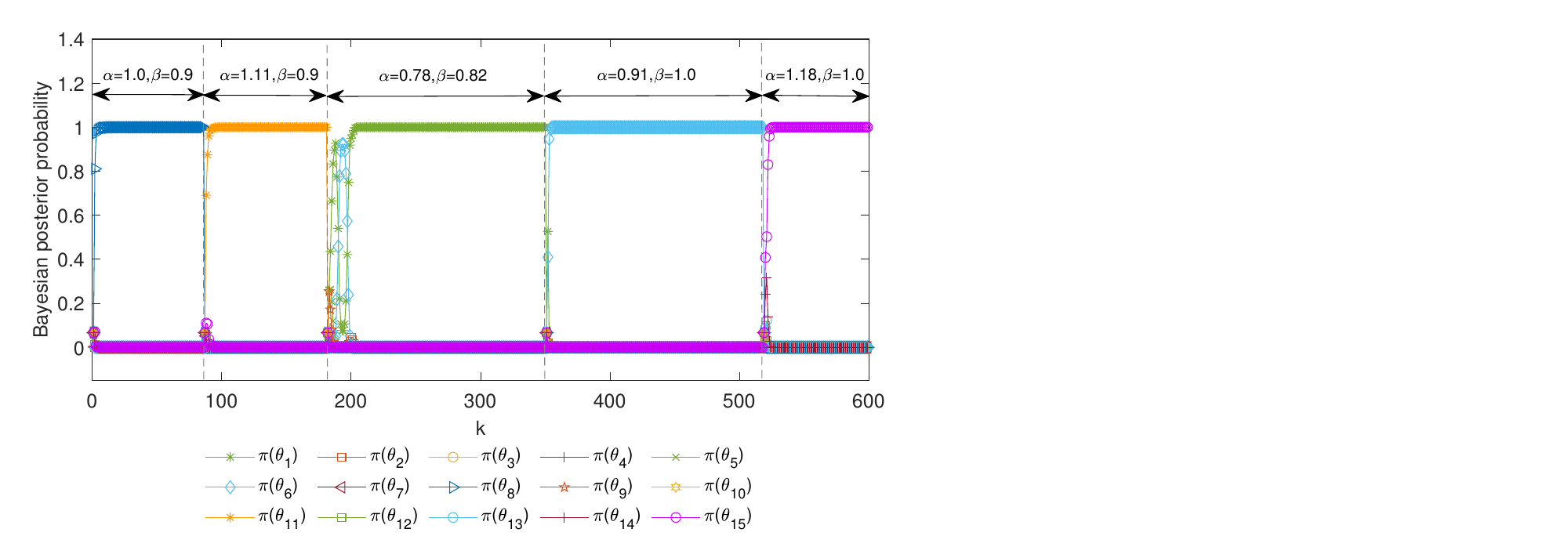}\\
	\caption{The convergence of Bayesian posterior probability for disturbance candidates over time in case 1.}\label{Fig_case1_1}
\end{figure}

\begin{figure}[htbp]
	\centering
	\includegraphics[width=0.50\textwidth]{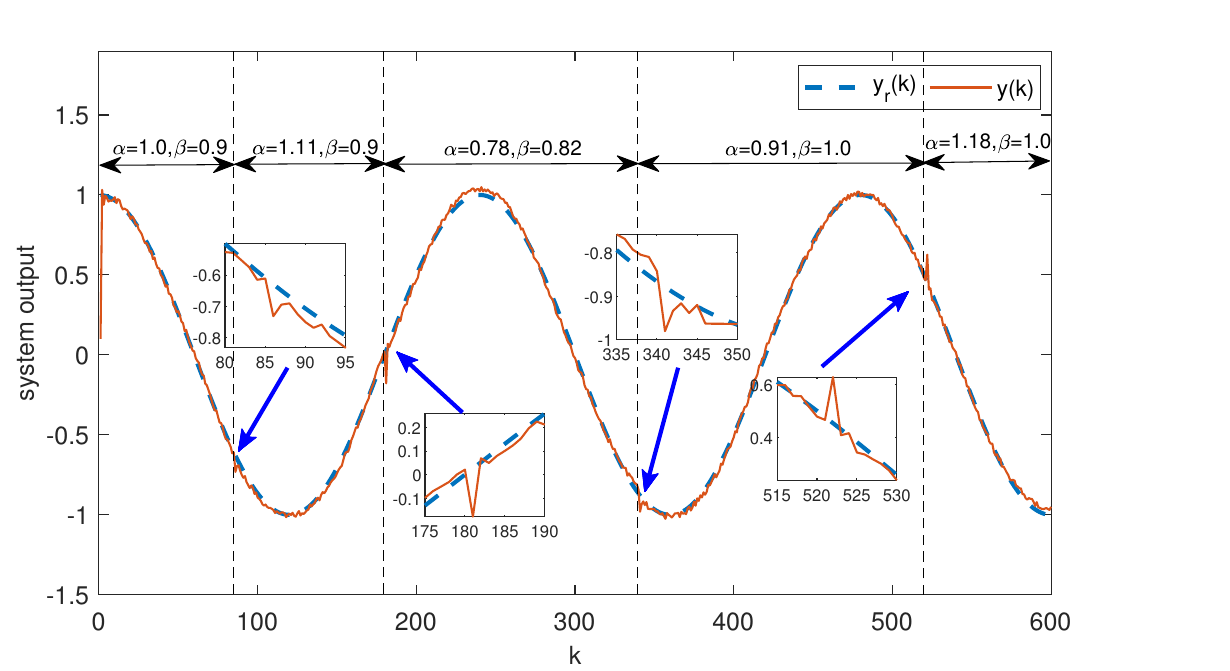}\\
	\caption{The system output tracking performance by our approach in case 1. $y_r(k)$ and $y(k)$ represent the reference trajectory and the system output, respectively.}\label{Fig_case1_2}
\end{figure}

\begin{figure}[htbp]
	\centering
	\includegraphics[width=0.50\textwidth]{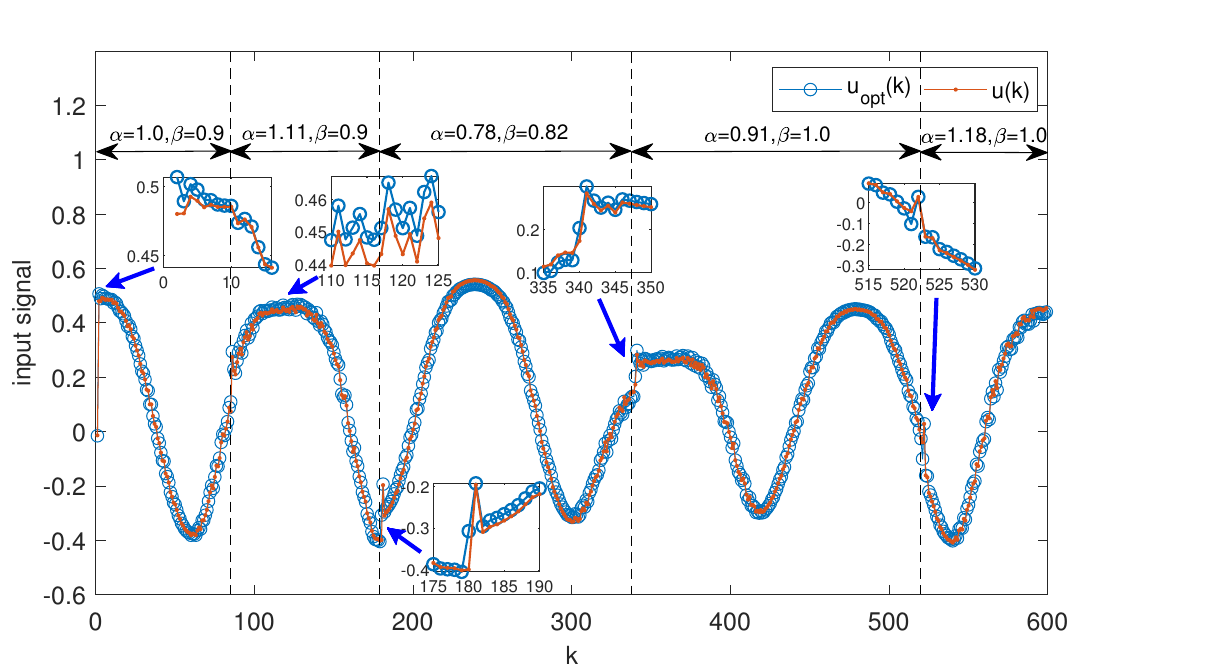}\\
	\caption{The control input signal calculated by our approach and the optimal control over time in case 1. $u_{opt}(k)$ and $u(k)$ represent the control signal by the optimal control and the proposed approach, respectively.}\label{Fig_case1_3}
\end{figure}

Fig.\ref{Fig_case1_1} shows the convergence for all Bayesian posterior probabilities with $\theta_t$-candidate ($t\in\{1,2,3,...,15\}$). We observe that during the process that the multiplicative disturbances are $\alpha=1.0$ and $\beta=0.9$, the posterior probability $\pi(\theta_8)$ converges from $1/15$ to $1$, and the other posterior probabilities $\pi(\theta_t)\ (t\neq 8)$ converge from $1/15$ to $0$. When the multiplicative disturbance $\alpha(k)$ varies from $1.0$ to $1.11$ at the $85$th iteration, all the Bayesian posterior probabilities are reset to $1/15$, indicating that this variation is successfully detected. After this reset, the active learning works again, and the posterior probability $\pi(\theta_8)$ converges to $0$, and the posterior probability $\pi(\theta_{11})$ converges to $1$ quickly. From the $180$th to the $340$th iteration, where the multiplicative disturbances change to $\alpha=0.78$ and $\beta=0.82$, the posterior probability $\pi(\theta_1)$ converges to $1$, meaning that the multiplicative disturbances $\alpha(k)$ and $\beta(k)$ are all located in the sub-interval $[0.75, 0.85]$, and are approximated to $0.8$. When the disturbances change to $\alpha=0.91$ and $\beta=1.0$ at the $340$th iteration, the posterior probability $\pi(\theta_6)$ converges to $1$, corresponding to that the $\alpha(k)$ and $\beta(k)$ are located in the sub-interval $[0.85, 0.95]$ and $[0.95, 1.05]$,and are approximated to $0.9$ and $1.0$ respectively. The posterior probability $\pi(\theta_{15})$ converges to $1$, and other probabilities converge to $0$ after the $520$th iteration as there is a new disturbance change at the $520$th iteration.

Fig.\ref{Fig_case1_2} depicts the system output tracking performance by our method. The blue dash line in this figure represents the reference trajectory, and the red line is the system output using the proposed control strategy. Notice that despite unknown disturbances, the output accurately tracks the reference trajectory under our approach, except for a short spike within $2-3$ iterations following the disturbances variation. This spike is unavoidable for the reason that the posterior probabilities need a short time to converge to the proper value when the disturbances change. %This is due to the unavailability of any exact information about the variation of the disturbances.

Fig.~\ref{Fig_case1_3} shows the control input signals calculated by our approach compared to that generated by the optimal control law. The optimal control law is derived with the system dynamics and disturbances completely known to the controller. Thus, compared with the proposed control law, the optimal control is a much more accurate strategy and free from the uncertainties induced by disturbances and can be regarded as the ideal control method to benchmark our proposed approach. Fig.~\ref{Fig_case1_3} shows that, as expected, our control input value converges to that of the optimal control law from the $1$st iteration up to the $85$th iteration, with the disturbance learned approximating the true value. The sub-figures in Fig.~\ref{Fig_case1_3} show the deviation between our control law and the optimal control law around the $1$st, $85$th, $180$th, $340$th, and $520$th iterations. Those derivations happen where there is a change in the disturbances. We also note that this deviation can deteriorate because the approximated disturbances are not the true values. I.e., there can be a small error between the true and the estimated disturbance. For example, from the $180$th to $340$th iteration, the true value of the disturbances are $\alpha=0.78$ and $\beta=0.82$, while the learned approximated values are $\hat{\alpha}=0.8$ and $\hat{\beta}=0.8$, respectively. The deviation of our approach from the optimal control could be reduced by dividing the bounded disturbance into more sub-intervals to obtain a more accurate disturbance approximation, which we will further discussed in Case 3 in Section~\ref{section4.3}.

\subsection{Case 2: simulations where the nonlinear system is disturbed by additive disturbance and the reference trajectory is a square wave}\label{section4.2}

Consider the nonlinear system described in equation~(\ref{case1}). Assume that the additive disturbance $\gamma(k)$ is already estimated offline to be in the bounded interval $[-1.45, 0.55]$, and the multiplicative disturbances are $\alpha(k)=1$ and $\beta(k)=1$. The system noise $e(k)$ follows Gaussian distribution $N(0,0.0025)$. The reference trajectory over time is set as a square wave
\begin{equation} \label{squarewave}
	\begin{aligned}
		y_r(k) = \left\{ \begin{array}{rcl}
			1 & \mbox{for}& 1 \leq k < 150\\
			-1 & \mbox{for} & 150 \leq k < 300\\
			1 & \mbox{for} & 300 \leq k < 450\\
			-1 & \mbox{for} & 450 \leq k \leq 600\\
		\end{array}\right.
	\end{aligned}
\end{equation}
The variation of the additive disturbance is
\begin{equation} \label{case2gamma}
	\begin{aligned}
		\gamma(k) = \left\{ \begin{array}{rcl}
			-0.2 & \mbox{for}& 1 \leq k < 200\\
			0.32 & \mbox{for} & 200 \leq k < 400\\
			-0.85 & \mbox{for} & 400 \leq k < 500\\
			0.5 & \mbox{for} & 500 \leq k \leq 600\\
		\end{array}\right.
	\end{aligned}
\end{equation}

We set the admissible additive disturbance approximation error as $\varepsilon_{\gamma}=0.1$, and the sub-interval number is $s_{\gamma}=20$. The corresponding disturbance approximation finite set is $\Omega_{\gamma}=\{-1.4, -1.3, -1.2, -1.1, -1.0, -0.9, -0.8, -0.7, -0.6, -0.5, \\-0.4, -0.3, -0.2, -0.1, 0, 0.1, 0.2, 0.3, 0.4, 0.5\}$. The number of elements in the finite disturbance vector set $\Omega$ is $20$. The initial Bayesian posterior probabilities are set as $\pi(\theta_t|\mathfrak{I}^0)=1/20$ for $t=1,2,\cdots, 20$. The unlock criterion parameter is set as $\epsilon=0.5$. The initial disturbance estimation error covariance is set as $P(1)=I_{3\times3}$, and the dual-property control coefficient is set as $\lambda=0.9$.

\begin{figure}[htbp]
	\centering
	\includegraphics[width=0.8\textwidth]{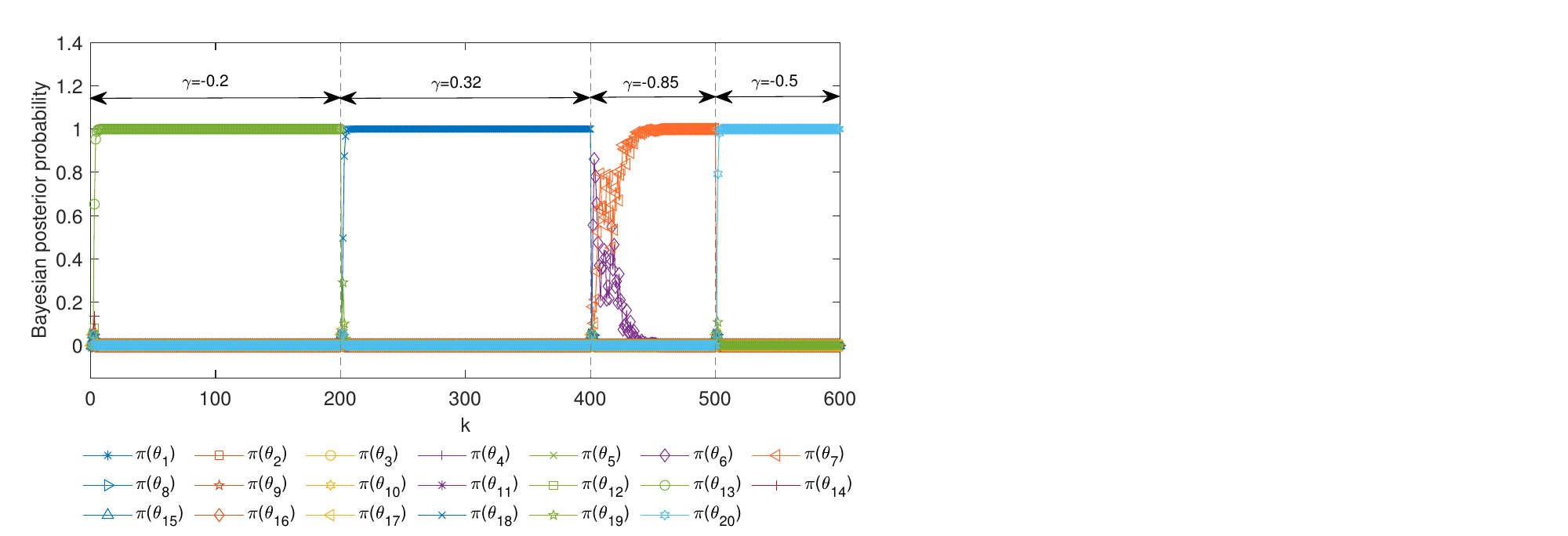}\\
	\caption{The convergence of Bayesian posterior probability for disturbance candidates over time in case 2.}\label{Fig_case2_1}
\end{figure}

Fig.\ref{Fig_case2_1} depicts the convergence of the Bayesian posterior probabilities for $\theta_t$-candidate. During the period $k\in[1,200]$ where the additive disturbance is $\gamma=-0.2$, the posterior probability $\pi(\theta_{13})$ converges from $1/20$ to $1$, and the posterior probability $\pi(\theta_{t})\ (t\neq 13)$ converges from $1/20$ to $0$, which means the disturbance $\gamma(k)$ is approximated to the right value $-0.2$. Notice that the additive disturbance $\gamma(k)$ varies from $-0.2$ to $0.32$ at the $200$th iteration and all the posterior probabilities are set to $\pi(\theta_i)=1/20$ once the variation is detected. The posterior probability $\pi(\theta_{18})$ converges to $1$ within five iterations after the disturbance variation at the $200$th iteration. This means that the additive disturbance $\gamma(k)$ is located in the sub-interval $[-0.35, -0.25]$, and we approximate the $\gamma(k)=0.32$ to $0.3$. From the $400$-$500$th iteration, the disturbance is $-0.85$, the edge point between the sub-interval $[-0.95, -0.85]$ and $[-0.85, -0.75]$. Therefore, the disturbance might converge to the sub-interval $[-0.95, -0.85]$ or $[-0.85, -0.75]$. As a result, the disturbance during this period may be approximated to $-0.9$ or $-0.8$. As expected, the curves for the posterior probability $\pi(\theta_{6})$ and $\pi(\theta_{7})$ fluctuate between $0$ and $1$ for within in a short time after the $400$th iteration, and finally $\pi(\theta_{6})$ converges to $0$ and $\pi(\theta_{7})$ converges to $1$. The disturbance $\gamma(k)$ changes to $-0.5$ at the $500$th iteration, and the  $\pi(\theta_{20})$ converges to $1$ after the disturbance variation.

\begin{figure}[htbp]
	\centering
	\includegraphics[width=0.75\textwidth]{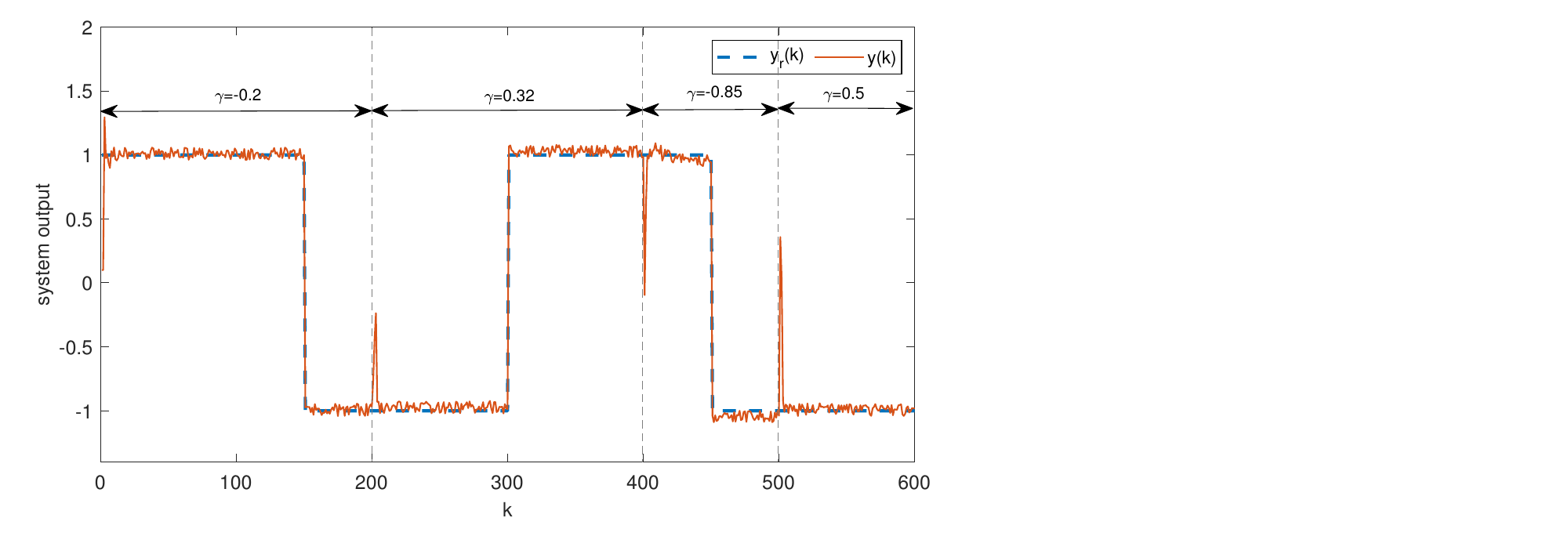}\\
	\caption{The system output tracking performance by our approach in case 2. $y_r(k)$ and $y(k)$ represent the reference trajectory and the system output, respectively.}\label{Fig_case2_2}
\end{figure}

Fig.\ref{Fig_case2_2} shows the system output tracking with a square wave trajectory by our approach. We can observe that under our approach, the system output accurately tracks the reference trajectory even for a step reference signal, except for a short spike within five iterations after the iterations where the disturbance changes. The spike is unavoidable since the posterior probability for the disturbance needs time to converge to the proper value after the disturbance varies.

\begin{figure}[htbp]
	\centering
	\includegraphics[width=0.75\textwidth]{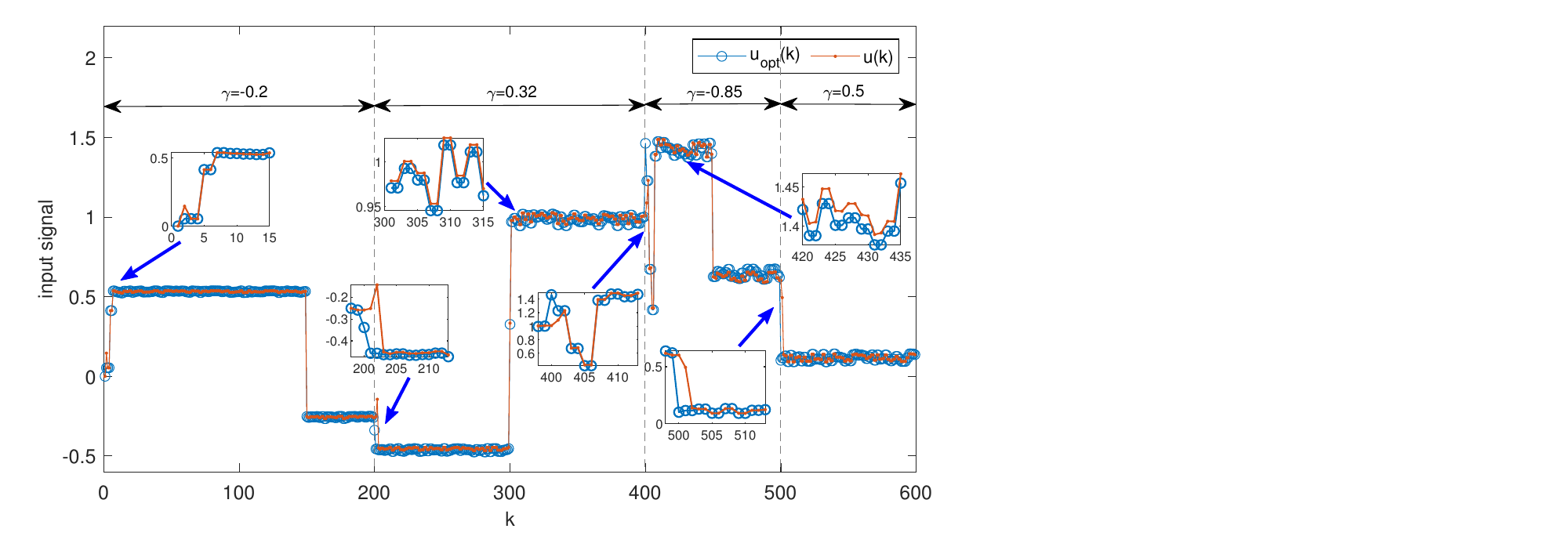}\\
	\caption{The control input signal calculated by our approach and the optimal control over time in case 2. $u_{opt}(k)$ and $u(k)$ represent the control signal by the optimal control and the proposed approach, respectively.}\label{Fig_case2_3}
\end{figure}

Fig.~\ref{Fig_case2_3} compares the control input signals by the proposed method with the ideal benchmark optimal control. It shows that the control signal by our method converges to the optimal control signal within $4$ iterations after the disturbance varies at the $1$st and $500$th iteration. The disturbance $\gamma(k)$ varies from $-0.2$ to $0.32$ at the $200$th iteration, and the control signal of our method converges to the control law with the disturbance approximation $\hat{\gamma}=0.3$ (see Fig.~\ref{Fig_case2_1}), which is close to the optimal control signal. We detail the deviation between our approach and the ideal benchmark in the sub-figure. When the additive disturbance changes to $\gamma(k)=-0.85$ at the $400$th iteration, as expected, the control signal fluctuates between the $\theta_6$-candidate and $\theta_7$-candidate control signal and gradually converges to the $\theta_6$-candidate control signal, which follows closely to the benchmark control signal by the optimal control.

\subsection{Case 3: Monte Carlo simulations regarding the admissible disturbance approximation error}\label{section4.3}
%\subsection{Case 3: Monte Carlo simulations regarding the admissible disturbance approximation error $\varepsilon_{\alpha}$, $\varepsilon_{\beta}$ and $\varepsilon_{\gamma}$}\label{section4.3}

We extend the simulations in case 1 and case 2 by running Monte Carlo simulations for the statistical experiment results of our approach. Simulations in this section keep the settings the same as those in case 1 and case 2, except that we use various admissible disturbance approximation errors $\varepsilon_{\alpha}$, $\varepsilon_{\beta}$, and $\varepsilon_{\gamma}$ for the disturbance $\alpha$, $\beta$, and $\gamma$, respectively. We vary those parameters to analyze how different resolutions of disturbance approximation can influence the performance of our approach.
%In this subsection, we compare the control performance for our proposed method with different values of parameters $\varepsilon_{\alpha}$, $\varepsilon_{\beta}$ and $\varepsilon_{\gamma}$. The simulation settings in this case are the same as that in case 1 and case 2.

In the Monte Carlo simulation for case 1, we test our proposed method with the parameters $\varepsilon_{\alpha}=\varepsilon_{\beta}=0.05$, $\varepsilon_{\alpha}=\varepsilon_{\beta}=0.075$, $\varepsilon_{\alpha}=\varepsilon_{\beta}=0.1$, and $\varepsilon_{\alpha}=\varepsilon_{\beta}=0.2$. When the parameters $\varepsilon_{\alpha}=\varepsilon_{\beta}=0.05$, the number of sub-intervals for $\alpha$ and is $\beta$ are $10$ and $6$, respectively. Then we have the corresponding disturbance approximation finite sets for $\alpha$ and $\beta$ as $\Omega_{\alpha}=\{0.775, 0.825, 0.875, 0.925, 0.975, 1.025, 1.075,\\ 1.125, 1.175, 1.225\}$ and $\Omega_{\beta}=\{0.775, 0.825, 0.875, 0.925, \\0.975, 1.025\}$. When the parameters $\varepsilon_{\alpha}=\varepsilon_{\beta}=0.075$, the number of sub-intervals for $\alpha$ and $\beta$ are $7$ and $4$, respectively. The corresponding disturbance approximation finite sets for $\alpha$ and $\beta$ are $\Omega_{\alpha}=\{0.7857, 0.8571,  0.9285, 0.9999, 1.0713,\\ 1.1427, 1.2142\}$ and $\Omega_{\beta}=\{0.7875, 0.8625, 0.9375, 1.0125\}$. When the parameters $\varepsilon_{\alpha}=\varepsilon_{\beta}=0.1$, the number of sub-intervals for $\alpha$ and $\beta$ are $5$ and $3$, respectively. The corresponding disturbance approximation finite sets for $\alpha$ and $\beta$ are $\Omega_{\alpha}=\{0.8, 0.9, 1.0, 1.1, 1.2\}$ and $\Omega_{\beta}=\{0.8, 0.9, 1.0\}$. When the parameters $\varepsilon_{\alpha}=\varepsilon_{\beta}=0.2$, the number of sub-intervals for $\alpha$ and $\beta$ are $3$ and $2$, respectively. The corresponding disturbance approximation finite sets for $\alpha$ and $\beta$ are $\Omega_{\alpha}=\{0.8334, 1.0001, 1.1667\}$ and $\Omega_{\beta}=\{0.825, 0.975\}$. Note that different values of $\varepsilon_{\alpha}$ and $\varepsilon_{\beta}$ lead to different sub-interval length, resulting in disturbance approximation with different resolutions.

We perform $100$ Monte Carlo runs for the control of the nonlinear system described in equation~\ref{case1}. The disturbances in this simulation are multiplicative, represented by $\alpha(k)$ and $\beta(k)$ located randomly in the bounded interval $[0.75,1.25]$ and $[0.75,1.05]$, respectively. To quantify the control performance for the $100$ Monte Carlo runs, we define the average performance index as
\begin{equation}
	\begin{aligned}
		J_M = \frac{1}{N_{MC}} \sum_{i=1}^{N_{MC}} \left( \frac{1}{N}\sum_{k=1}^{N} \sqrt{[y(k)-y_r(k)]^2}\right) _i
	\end{aligned}
\end{equation}
where $N_{MC}$ is the number of Monte Carlo runs. In this case $N_{MC}$ is $100$.

\begin{table}[H]
	\caption{Comparison of the average performance index $J_M$ and running time $t_{r}$ between the ideal benchmark optimal control and the proposed method with different $\varepsilon_{\alpha}$ and $\varepsilon_{\beta}$ (for the Monte Carlo simulation of case 1).}
	\centering
	\resizebox{4.5cm}{!}
	{\begin{tabular}{cccc}
			\toprule
			$\varepsilon_{\alpha}$, $\varepsilon_{\beta}$ &  $s_{\alpha}s_{\beta}$ & $J_M$ & $t_{r} (s)$\\
			\midrule
			0.05 & 60 & 0.0022 & 0.4110\\
			0.075 & 28 & 0.0027 & 0.1817\\
			0.1 & 15 & 0.0029 & 0.0954\\
			0.2 & 6 & 0.0042 & 0.0392\\
			\multicolumn{2}{c}{optimal control} & 0.0018 & 0.0036 \\
			\bottomrule
	\end{tabular}}
	\label{table2}
\end{table}

Table~\ref{table2} shows the results of Monte Carlo simulations for different $\varepsilon_{\alpha}$ and $\varepsilon_{\beta}$. The average performance index with $\varepsilon_{\alpha}=\varepsilon_{\beta}=0.05$ is the lowest among these different parameters, and the average performance index grows with the value of parameters $\varepsilon_{\alpha}$ and $\varepsilon_{\beta}$. It indicates that control performance with $\varepsilon_{\alpha}=\varepsilon_{\beta}=0.05$ is the best and closest to the ideal benchmark optimal control. Nevertheless, as anticipated, such an enhanced control performance under low $\varepsilon_{\alpha}$ and $\varepsilon_{\beta}$ is at the cost of additional computation. As shown in the table, the running time $t_r$ grows almost linearly with the number of disturbance candidates calculated by $s_{\alpha}s_{\beta}$.

In the Monte Carlo simulation for case 2, we use different admissible disturbance approximation error for the additive disturbance as $\varepsilon_{\gamma}=0.05$, $\varepsilon_{\gamma}=0.1$, $\varepsilon_{\gamma}=0.2$, and $\varepsilon_{\gamma}=0.4$. When the parameter $\varepsilon_{\gamma}=0.05$, the number of sub-intervals for $\gamma$ is $40$, and the corresponding disturbance approximation finite set is $\Omega_{\gamma} = \{-1.425, -1.375, -1.325,-1.275, -1.225,\\ -1.175, -1.125, -1.075, -1.025, -0.975, -0.925, -0.875,\\ -0.825, -0.775, -0.725, -0.675, -0.625, -0.575, -0.525,\\ -0.475, -0.425, -0.375, -0.325, -0.275, -0.225, -0.175,\\ -0.125, -0.075, -0.025, 0.025, 0.075, 0.125, 0.175, 0.225,\\ 0.275, 0.325, 0.375, 0.425, 0.475, 0.525\}$. When the parameter $\varepsilon_{\gamma}=0.1$, the number of sub-intervals for $\gamma$ is $20$, and the corresponding disturbance approximation finite set is $\Omega_{\gamma} = \{-1.4, -1.3, -1.2, -1.1, -1.0, -0.9, -0.8, -0.7, -0.6, -0.5,\\ -0.4, -0.3, -0.2, -0.1, 0, 0.1, 0.2, 0.3, 0.4, 0.5 \}$. When the parameter $\varepsilon_{\gamma}=0.2$, the number of sub-intervals for $\gamma$ is $10$, and the corresponding disturbance approximation finite set is $\Omega_{\gamma} = \{ -1.35, -1.15, -0.95, -0.75, -0.55, -0.35, -0.15,\\ 0.05, 0.25, 0.45 \}$. When the parameter $\varepsilon_{\gamma}=0.4$, the number of sub-intervals for $\gamma$ is $5$, and the corresponding disturbance approximation finite set is $\Omega = \{ -1.25, -0.85, -0.45, -0.05, \\0.35\}$. Note that different values of $\varepsilon_\gamma$ lead to a disturbance approximation finite set of different sizes, meaning that the disturbance is approximated to different levels.

We perform $100$ Monte Carlo runs to control the nonlinear system described in equation~\ref{case1} corrupted by additive disturbance. The additive disturbance $\gamma(k)$ is random within the bounded interval $[-1.45,0.55]$. Table~\ref{table3} shows the results of Monte Carlo simulations for different $\varepsilon_{\gamma}$. As anticipated, the average performance index with $\varepsilon_{\gamma}=0.05$ is the lowest among these different parameters and closest to that of the optimal control, and the average performance index grows with the parameter $\varepsilon_{\gamma}$, indicating a better control performance for a lower $\varepsilon_{\gamma}$. Similar to the Monte Carlo simulation for case 1, we can also conclude that lower $\varepsilon_{\gamma}$ leads to more disturbance candidates that incur heavier computation, which is indicated by the calculation time for the same control iterations in the table. 

\begin{table}[H]
	\caption{Comparison of the average performance index $J_M$ and running time $t_{r}$ between the ideal benchmark optimal control and our proposed method with different $\varepsilon_{\gamma}$ (for the Monte Carlo simulation of case 2).}
	\centering
	\resizebox{4.5cm}{!}
	{\begin{tabular}{cccc}
			\toprule
			$\varepsilon_{\gamma}$ &  $s_{\gamma}$ & $J_M$ & $t_{r} (s)$ \\
			\midrule
			0.05 & 40 & 0.0030 & 0.3021\\
			0.1 & 20 & 0.0033 & 0.1472\\
			0.2 & 10 & 0.0037 & 0.0732\\
			0.4 & 5 & 0.0054 & 0.0377\\
			\multicolumn{2}{c}{optimal control}  &  0.0024 & 0.0036\\
			\bottomrule
	\end{tabular}}
	\label{table3}
\end{table}

\subsection{Case 4: evaluation of our approach in the speed control of high-speed train} \label{section4.4}
%\subsection{Case 4: Extensive comparison simulations for the speed control of high-speed train} \label{section4.4}

This evaluation considers the speed control of CRH3 high-speed train. We evaluate our approach in this case to check its potential applicability where the dynamics of a real nonlinear system is unknown to the controller design and corrupted by disturbances varying over time. We implement our approach along with two other control methods for performance comparison, i.e., the ideal benchmark optimal control and a latest model-free adaptive control method. 

The dynamics of a CRH3 high-speed train can be described as~\cite{chou2007modelling}:
\begin{equation}\label{trainmodel}
	\begin{aligned}
		v(k+1) = v(k)+\xi T \{ f(k)-W[v(k)]\}\\
		%	\end{aligned}
	%\end{equation}
	%\begin{equation}
	%	\begin{aligned}
		W[v(k)] = c_r + c_mv(k) + c_av^2(k)
	\end{aligned}
\end{equation}
where $v(k)$ is the speed of the high-speed train, $f(k)$ is the traction/braking force, $\xi$ is the acceleration coefficient, $T$ is the sampling interval, $W[v(k)]$ is the general resistance, which includes the rolling resistance $c_r$, the mechanical resistance $c_mv(k)$, and the aerodynamic drag $c_av^2(k)$. The speed control for high-speed train aims to keep the train at an objective running speed by adjusting the traction/braking force. However, since the train could be disturbed by the abnormal train tracks condition, the changing weather and climate, and mechanical wear and aging, the speed control of the train can be negatively affected. Beyond that, the coefficients $\xi$, $c_r$, $c_m$, and $c_a$ are not constant and would be disturbed by the complex outer environmental and inner components conditions.

In the implementation of our approach for this case, we first reconstruct the train model in equation~(\ref{trainmodel}) to fit the structure of our proposed specialized neural network (SNN) as:
\begin{equation} \label{hstrain}
	\begin{aligned}
		v(k+1) = & \alpha(k) \{ v(k)-\xi T [c_r + c_mv(k) + c_av^2(k)] \} \\
		&+ \beta(k) \xi T f(k) + \gamma(k) + e(k+1)
	\end{aligned}
\end{equation}
where we include the multiplicative disturbance $\alpha(k)$ and $\beta(k)$ and the additive disturbance $\gamma(k)$. Since the train dynamics is supposed to be unknown to the control design, the SNN uses two nonlinear functions to approximate the nonlinear dynamics in equation~(\ref{hstrain}) that $\hat{f}[x(k)] \approx v(k)-\xi T [c_r + c_mv(k) + c_av^2(k)]$ and $\hat{g}[x(k)] \approx \xi T$. The parameters for the function $\hat{f}[x(k)]$ and $\hat{g}[x(k)]$ can be obtained through offline learning from historical measurements of the train.

For the multiplicative and additive disturbance, we can estimate empirically from the system measurements that $\alpha(k)$, $\beta(k)$, are $\gamma(k)$ are bounded within the interval $[0.725,1.075]$, $[0.8,1.2]$, and $[-13,2]$, respectively. We set the random noise $e(k)$ following the Gaussian distribution $N(0,1)$. We use the following parameters of a high-speed train to conduct this simulation. The parameters of a CRH3 train for equation~(\ref{hstrain}) are $\xi=0.06$, $c_r=0.1$, $c_m=0.0064$ and $c_a=0.000115$, and the sampling interval is $T=0.1$. We set the disturbances over time to be:
\begin{equation} \label{trainmulalpha}
	\begin{aligned}
		\alpha(k) = \left\{ \begin{array}{rcl}
			1.0 & \mbox{for}& 1 \leq k < 100\\
			1.05 & \mbox{for}& 100 \leq k < 250\\
			0.95 & \mbox{for} & 250 \leq k < 350\\
			0.92 & \mbox{for} & 350 \leq k < 500\\
			0.95 & \mbox{for} & 500 \leq k \leq 600\\
		\end{array}\right.
	\end{aligned}
\end{equation}
\begin{equation} \label{trainmulbeta}
	\begin{aligned}
		\beta(k) = \left\{ \begin{array}{rcl}
			1.0 & \mbox{for}& 1 \leq k < 100\\
			0.9 & \mbox{for}& 100 \leq k < 350\\
			1.0 & \mbox{for} & 350\leq k < 500\\
			1.15 & \mbox{for} & 500 \leq k \leq 600\\
		\end{array}\right.
	\end{aligned}
\end{equation}
\begin{equation} \label{trainmulgamma}
	\begin{aligned}
		\gamma(k) = \left\{ \begin{array}{rcl}
			0.0 & \mbox{for}& 1 \leq k < 350\\
			-12.5 & \mbox{for} & 350\leq k < 500\\
			-2.25 & \mbox{for} & 500\leq k \leq 600\\
		\end{array}\right.
	\end{aligned}
\end{equation}
where $k$ is the control iteration index. The reference speed for the control is set to vary over time as $270+50(1+\exp(-2k))$.

We set the parameters in our proposed method as $\varepsilon_{\alpha}=0.05$, $\varepsilon_{\beta}=0.4$ and $\varepsilon_{\alpha}=1$. The corresponding number of sub-interval for the disturbance candidates is $s_{\alpha}=7$, $s_{\beta}=1$ and $s_{\gamma}=15$. The disturbance approximation finite sets for $\alpha$, $\beta$, and $\gamma$ are $\Omega_{\alpha}=\{0.75, 0.80, 0.85, 0.90, 0.95, 1.00, 1.05\}$, $\Omega_{\beta}=\{1.00\}$, and $\Omega_{\gamma} = \{ -12.5, -11.5, -10.5, -9.5, -8.5, -7.5, -6.5, -5.5, \\-4.5, -3.5, -2.5, -1.5, -0.5, 0.5, 1.5 \}$. The number of the candidate disturbance vector $\theta_t$ is $105$, and the initial Bayesian posterior probabilities are set as $\pi(\theta_t|\mathfrak{I}^0)=1/105$ for $t=1,2,\cdots, 105$. The unlock criterion parameter is set as $\epsilon=1.5$. The initial disturbance estimation error covariance is set as $P(1)=I_{3\times3}$, and the dual-property control coefficient is set as $\lambda=0.9$.

\begin{figure}[htbp]
	\centering
	\includegraphics[width=0.4\textwidth]{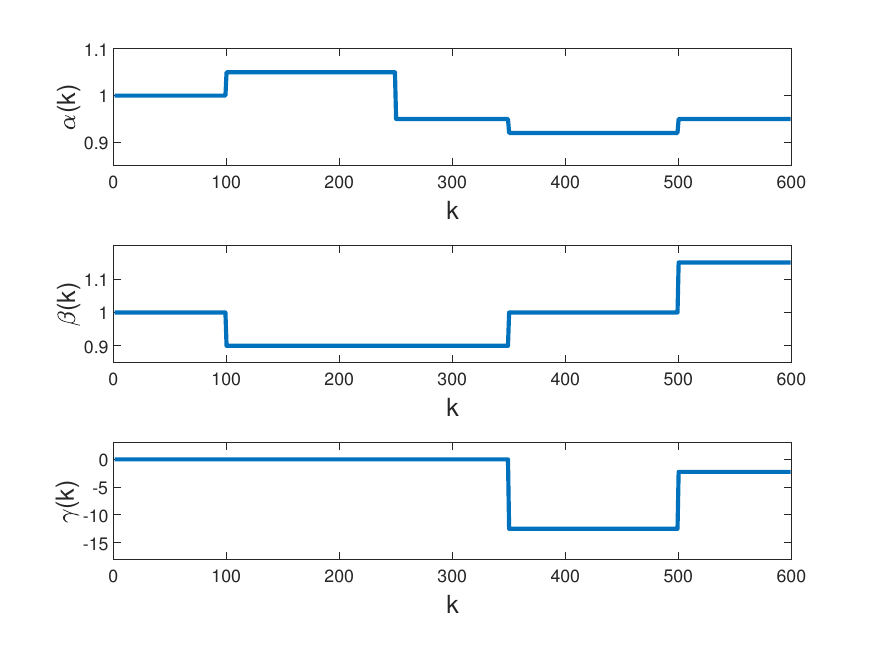}\\
	\caption{The multiplicative disturbances $\alpha(k)$ and $\beta(k)$, and the additive disturbance $\gamma(k)$ over time for the evaluation in case 4.}\label{Fig_case4_1}
\end{figure}

\begin{figure}[htbp]
	\centering
	\includegraphics[width=0.4\textwidth]{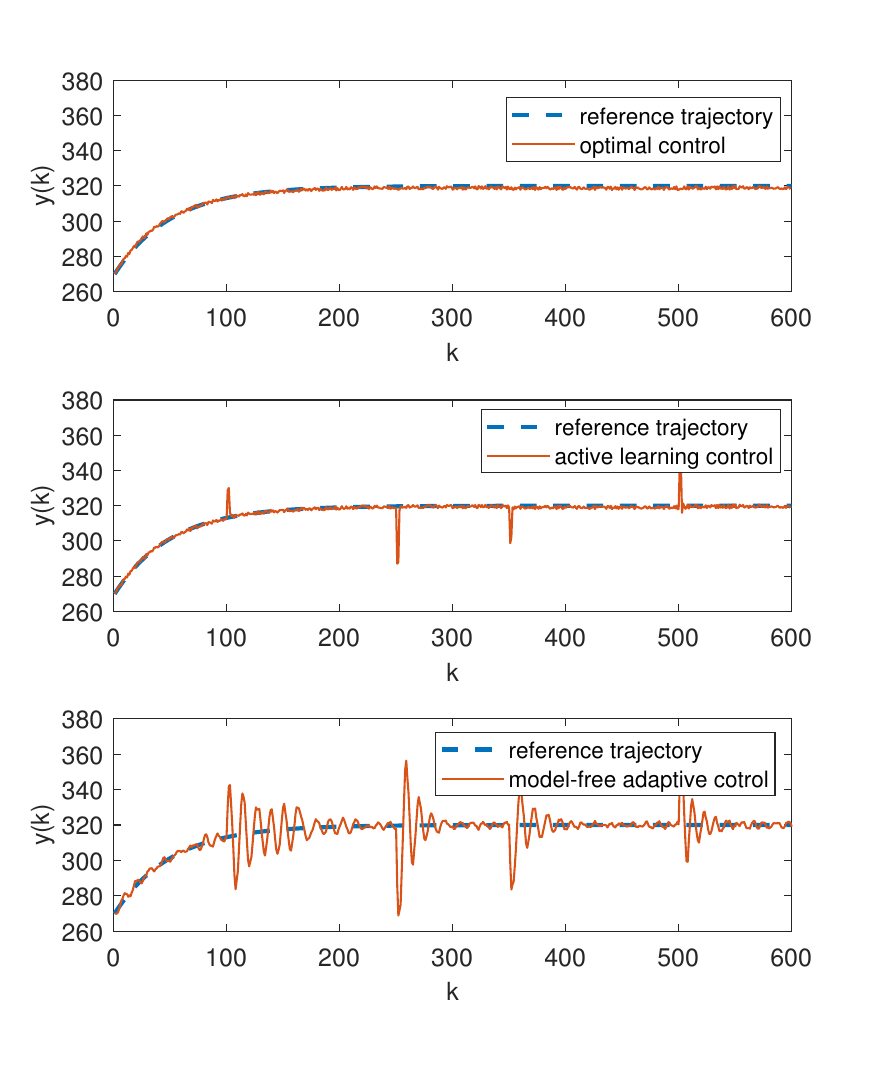}\\
	\caption{The system output of the high-speed train by the optimal control, the proposed approach, and the model-free adaptive control}\label{Fig_case4_2}
\end{figure}

\begin{figure}[htbp]
	\centering
	\includegraphics[width=0.4\textwidth]{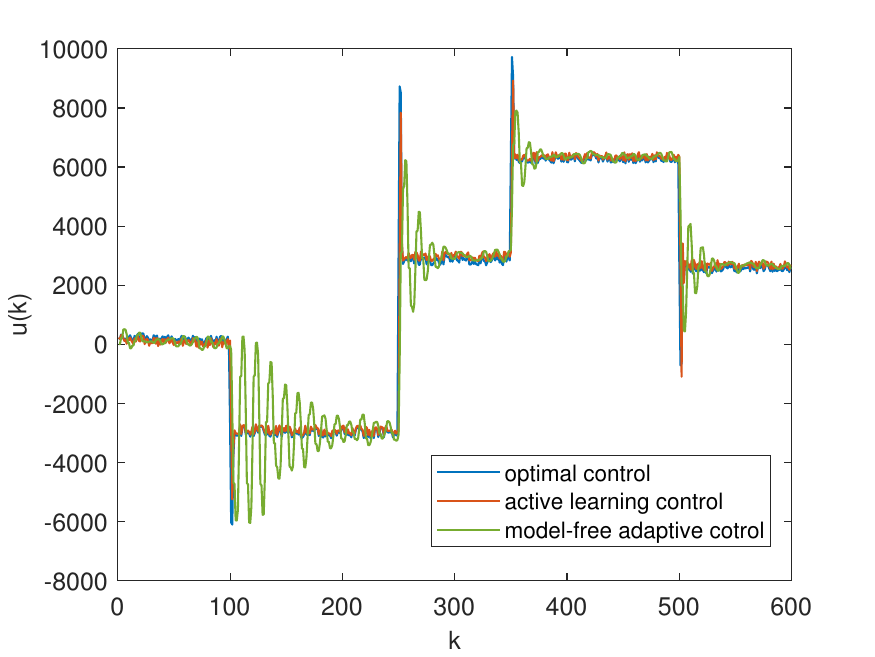}\\
	\caption{The system control input signals by the optimal control, the proposed approach, and the model-free adaptive control }\label{Fig_case4_3}
\end{figure}

\begin{figure}[htbp]
	\centering
	\includegraphics[width=0.45\textwidth]{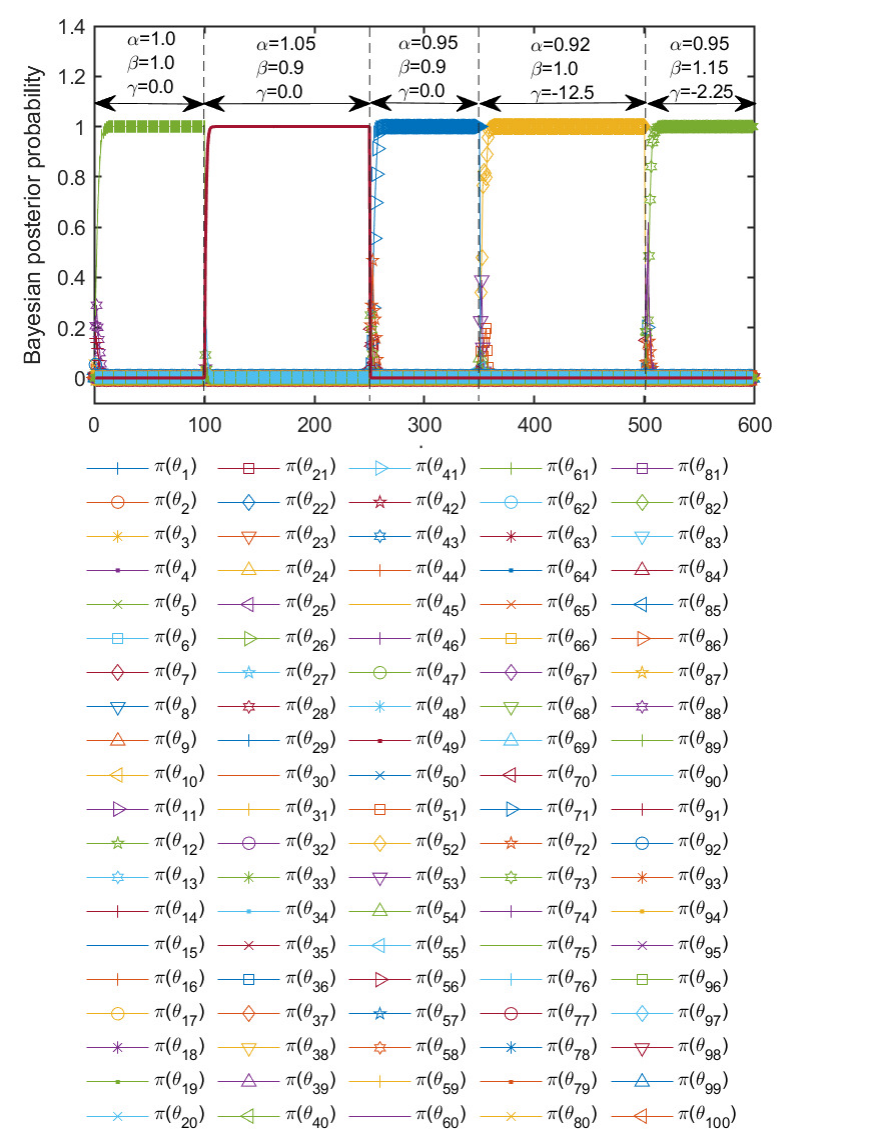}\\
	\caption{The convergence of Bayesian posterior probability for the disturbance candidates during the active learning of our approach.}\label{Fig_case4_4}
\end{figure}

\begin{figure}[htbp]
	\centering
	\includegraphics[width=0.4\textwidth]{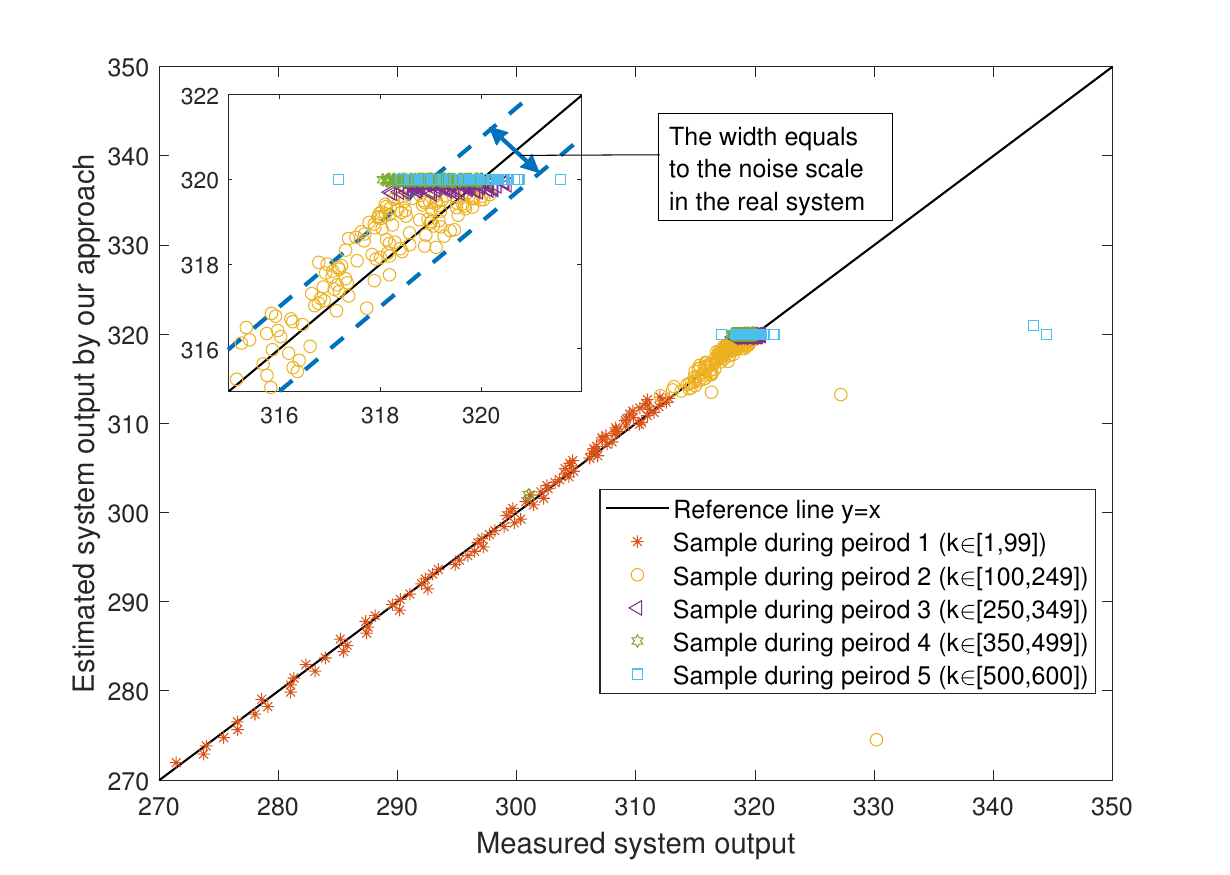}\\
	\caption{The comparison of measured system output and the estimated system output. $k$ in the legend is the control iteration index.}\label{Fig_case4_5}
\end{figure}

Fig.~\ref{Fig_case4_1} depicts how the disturbances $\alpha(k)$, $\beta(k)$, and $\gamma(k)$ vary over time. We run this evaluation in $600$ iterations, and Fig.~\ref{Fig_case4_2} shows the system output $y(k)$ under the proposed method compared with the ideal benchmark optimal control and the model-free adaptive control~\cite{6093751}. We can observe that the optimal control provides a perfect output tracking performance since the system disturbances are known to the control law derivation. The proposed approach also shows a good tracking control performance except for a short spike within $3-4$ iterations following a change in the disturbances. In comparison, the model-free control brings large overshoots and a long settling time (ranging from ca. 50 to 100 iterations) following the disturbances variation. Fig.~\ref{Fig_case4_3} compares the control input signal by the ideal benchmark optimal control, the proposed approach, and the model-free adaptive control, where it manifests that our approach provides the control signal much closer to the ideal benchmark and is more resilient to disturbances than the model-free adaptive control. Fig.~\ref{Fig_case4_4} shows the convergence of Bayesian posterior probability $\pi(\theta_t)$ during our approach's active learning. We observe that the posterior probability for new disturbance candidates converges quickly after a change in the disturbance. We use the disturbances learned in the active learning to estimate the system output $\hat{y}(k)$, and compare the estimated system output with the ground truth in Fig.~\ref{Fig_case4_5}. The comparison demonstrates that most of the estimated system output is close to the measured output. We observe that the divergences between the estimated system output and the ground truth are mostly bounded within the interval $[-1,1]$, which equals the scale of the Gaussian noise, as shown in equation~(\ref{hstrain}). Such a bounded system output error indicates the efficiency of the active learning where the Bayesian posterior probability of the disturbance candidates converges to the right value. Moreover, it also manifests that the disturbed nonlinear system is accurately approximated during the control law derivation.

\section{Conclusion} \label{section5}

This work studies the control of unknown nonlinear systems corrupted by disturbance. We propose an anti-disturbance dual control for such systems, where we learn the disturbance in real-time to reduce its influence on the control law derivation. In the proposed solution, we consider the influence of both multiplicative and additive disturbance, and we design a specialized neural network (SNN) to approximate the unknown system dynamics and disturbances numerically. Furthermore, while our approach executes the disturbance learning in the SNN and the system output tracking simultaneously and iteratively, we design the active learning for the disturbance to decouple from the output tracking, which could enhance the control robustness in case of varying and abrupt disturbances in the system. Numerical simulations and evaluations on the speed control of high-speed train demonstrate the efficiency of our disturbance learning strategy, the resilience of our approach against disturbances, and robust control performance.

In future work, we plan to extend our approach to multiple-input-multiple-output systems that better fit practical cases. Meanwhile, we anticipate to apply the proposed approach to the control of practical systems, e.g., the automatic control of remotely operated vehicle or underwater machines, which may come up with complex disturbances exploring uncertain environments like an ocean or a river.

\bibliographystyle{unsrt}
\bibliography{bibfile}

\end{document}